\documentclass[preprint,12pt,authoryear]{elsarticle}

\usepackage{amssymb}
\usepackage{amsmath} 
\usepackage{graphicx}
\usepackage{hyperref}
\usepackage{booktabs}
\usepackage{listings}
\usepackage{xcolor}
\usepackage{ragged2e}
\usepackage{float}
\usepackage{algorithm} 
\usepackage{algpseudocode} 
\usepackage{caption}
\usepackage{setspace}
\usepackage[utf8]{inputenc}
\usepackage[modulo]{lineno}
\usepackage{tikz}
\usepackage{listings}
\usepackage{xcolor}
\usepackage{subcaption}
\usepackage{booktabs, multirow}
\usetikzlibrary{positioning}

\journal{}

\begin{document}

\begin{frontmatter}



\title{Elasto-plastic cell-based smoothed finite element method solving geotechnical problems}


\author[inst1]{Yang Yang}
\author[inst2]{Mingjiao Yan}
\author[inst1]{Zongliang Zhang}
\author[inst2]{Miao Zhang}
\author[inst1]{Feidong Zheng}
\author[inst1]{Dong Pan}
\author[inst1]{Xiaozi Lin}
\affiliation[inst1]{organization={PowerChina Kunming Engineering Corporation Limited},
            city={Kunming},
            postode={650051}, 
            state={Yunnan}, 
            country={China}}
\affiliation[inst2]{organization={College of Water Conservancy and Hydropower Engineering, Hohai University},
            city={Nanjing},
            postcode={210098}, 
            state={Jiangsu},
            country={China}}

\begin{abstract}
This work develops an elasto-plastic cell-based smoothed finite element method (CSFEM) for geotechnical analysis. The formulation incorporates a smoothed strain field into the standard elasto-plastic framework based on the Mohr–Coulomb criterion and is implemented in ABAQUS through a user-defined element (UEL). A UEL–UMAT data-transfer strategy is introduced to enable post-processing of stress and strain in ABAQUS. The method is assessed using several benchmark problems, including three classical examples, a tunnel excavation, and a slope stability analysis. The results show that the CSFEM achieves accuracy comparable to or slightly better than the conventional FEM and matches analytical or reference solutions for the examined cases. These findings indicate that the proposed CSFEM provides a practical and robust alternative for routine elasto-plastic analyses in geotechnical engineering.
\end{abstract}


\begin{highlights}
    \item Elasto-plastic CSFEM for geotechnical analysis is developed.
    \item Implement the method in ABAQUS via a UEL with a UEL–UMAT transfer scheme.
    \item Benchmarks show accuracy comparable to FEM with smoother contours.
\end{highlights}


\end{frontmatter}


\section{Introduction}
\label{sec:Introduction}
The reliable prediction of the mechanical response of geomaterials under complex loading conditions remains a fundamental challenge in geotechnical engineering. Geomaterials, such as soil and rock, often exhibit nonlinear behaviors including plasticity, strain-softening, and path-dependent responses~\citep{chu1992strain,kok2009case}. These features require robust numerical methods to capture the intricate interplay between elasticity and plasticity in large-scale simulations of foundations~\citep{ray2021application}, slopes~\citep{kardani2021improved}, tunnels~\citep{jin2022effect}, embankments\citep{xia2024stability}, and underground structures\citep{hemeda2022geotechnical}. Among various computational approaches, the finite element method (FEM) has become the most widely adopted due to its flexibility in handling irregular geometries and boundary conditions~\citep{bottero1980finite,liu2015slope,wei20253d}.

However, FEM formulations present several drawbacks in the elasto-plastic analysis of geotechnical problems. FEM is known to overestimate the stiffness matrix, leading to an overall response that is too stiff and an underestimation of the internal strain energy~\citep{fraeijs1965displacement}. As a result, the computed displacements are generally smaller than the actual values. In addition, the required mapping and coordinate transformation restrict element shapes: for a four-node isoparametric element, all interior angles must remain below 180$^{\circ}$ and the Jacobian determinant must remain positive~\citep{liuFiniteElementMethod2003,liuTheoreticalAspectsSmoothed2007,liuSmoothedFiniteElement2007}. These constraints increase computational cost and significantly limit FEM performance in geotechnical problems involving large deformation or severe mesh distortion.

By introducing the strain-smoothing technique first proposed by Liu et al.~\citep{liuSmoothedFiniteElement2007}, the smoothed finite element method (SFEM)~\citep{liu2016smoothed} has developed into a promising framework for overcoming the intrinsic limitations of classical FEM, thereby drawing substantial attention in recent years. The SFEM delivers higher accuracy and faster energy convergence than the conventional FEM~\citep{daiNsidedPolygonalSmoothed2007,yanFastCellSmoothedFinite2025}. Specifically, its freedom from mapping or coordinate transformation allows elements of arbitrary shapes, making it particularly beneficial for geotechnical analyses where large deformations can cause severe mesh distortion.

Several SFEM variants have been proposed, including face-based~\citep{huang2023face,jiang2025strength}, node-based~\citep{sun2022stable,lyu2024implicit}, edge-based~\citep{he2022edgeon,he2022edgean}, and cell-based formulations~\citep{surendran2021cell,liu2022cell,sun2025cell}.Among these, the cell-based SFEM (CSFEM) provides a particularly straightforward framework that remains close to the standard FEM procedure. In the CSFEM, each element is subdivided into smoothing sub-cells where the strain field is averaged, and the stiffness contribution is evaluated through surface integration based on the divergence theorem~\citep{lee2015three}.

At present, the CSFEM has been successfully applied to several physical field problems: Surendran et al.\citep{surendran2021cell} presents a CSFEM using polygonal elements to model interfacial cracks with non-matching grids, enabling accurate fracture analysis while reducing degrees of freedom and relaxing meshing constraints. Yan et al. \citep{yanFastCellSmoothedFinite2025} present a fast CSFEM to solve static and dynamic problems.
Lee et al. \citep{lee2015three} a three-dimensional CSFEM for application to elasto-plastic analysis for the solid mechanics. Peng et al. \citep{peng2021phase} developed a novel phase-field model for brittle fracture using CSFEM with spectral decomposition to capture nonlinear stress and elastic responses. These studies demonstrate the effectiveness and versatility of CSFEM, yet its application in geotechnical engineering simulations remains limited.

In this work, we propose an elasto-plastic CSFEM tailored for geotechnical applications. The method incorporates strain smoothing into the standard FEM formulation, replacing the compatible strain with a smoothed strain field evaluated over subcells and coupling it with a Mohr--Coulomb elasto-plastic model. This leads to a softer and more stable stiffness matrix, reduces mesh-distortion sensitivity, and improves stress prediction in nonlinear analyses. The remainder of this paper is organized as follows. Section~\ref{sec:Cell-based SFEM formulation for elasto-plasticity} presents the theoretical formulation of the elasto-plastic CSFEM. Section~\ref{Elasto-plastic constitutive formulation} introduces the Mohr--Coulomb constitutive model. Section~\ref{Construction of CSFEM elements} describes the construction of smoothing subcells and the smoothed strain–displacement matrix. Section~\ref{Implementation} outlines the ABAQUS UEL implementation and the UEL–UMAT data-transfer scheme. Section~\ref{Numerical examples} provides benchmark examples to assess accuracy and convergence. Finally, Section~\ref{Conclusions} summarizes the main findings.

\section{CSFEM formulation for elasto-plasticity}
\label{sec:Cell-based SFEM formulation for elasto-plasticity}
In this study, we examine a deformable body that occupies the domain $\Omega$. This body responds to internal body forces $\textbf{b}$, externally applied tractions $\overline{\textbf{t}}$ on the boundary $\Gamma_t$, and displacement boundary conditions $\mathbf{u}=\overline{\mathbf{u}}$ on $\Gamma_u$. The governing equation and boundary conditions are expressed as follows:

\begin{equation}
\begin{array}{rlc}
\nabla \boldsymbol{\sigma}+\mathbf{b}=0 & \text { in } & \Omega \\
\mathbf{u}=\overline{\mathbf{u}} & \text { in } & \Gamma_u, \\
\boldsymbol{\sigma} \mathbf{n}=\overline{\mathbf{t}} & \text { in } & \Gamma_t
\end{array}
\end{equation}
where $\nabla$ denotes the differential operator, $\boldsymbol{\sigma}$ represents the Cauchy stress tensor.

The relationship of stress $\boldsymbol{\sigma}$ and strain $\boldsymbol{\varepsilon}$ can be expressed by
\begin{equation}
    \boldsymbol{\sigma}=\boldsymbol{\mathrm{D}}\boldsymbol{\varepsilon},
\end{equation}
where $\boldsymbol{\mathrm{D}}$ is the constitutive matrix.

Mathematically expressing the virtual displacement and compatible strain within each element of the spatial discretization, we obtain the subsequent equations:
\begin{equation}
\delta \mathbf{u}^h=\sum_i \mathbf{N}_i \delta \mathbf{d}_i, \quad \mathbf{u}^h=\sum_i \mathbf{N}_i \mathbf{d}_i,
\end{equation}

\begin{equation}
\delta \boldsymbol{\varepsilon}^h=\sum_i \mathbf{B}_i \delta \mathbf{d}_i, \quad \boldsymbol{\varepsilon}^h=\sum_i \mathbf{B}_i \mathbf{d}_i,
\end{equation}
where $\mathbf{d}_i=\left[\begin{array}{ll}u_i & v_i\end{array}\right]^{\mathrm{T}}$  is the nodal displacement vector, and $\mathbf{B}_i$ is the strain-displacement matrix as follows:
\begin{equation}
\mathbf{B}_i=\left[\begin{array}{cc}
N_{i x} & 0 \\
0 & N_{i y} \\
N_{i y} & N_{i x}
\end{array}\right],
\end{equation}
where $N_{i x}$ and $N_{i y}$ are outer normal derivatives in relation to $x$ and $y$. The process of energy assembly yields the following results:
\begin{equation}
\int_{\Omega} \delta \mathbf{d}^{\mathrm{T}} \mathbf{B}^{\mathrm{T}} \boldsymbol{\sigma} \mathrm{d} \Omega-\int_{\Omega} \delta \mathbf{d}^{\mathrm{T}} \mathbf{N}^{\mathrm{T}}\mathbf{b}\mathrm{d} \Omega-\int_{\Gamma_t} \delta \mathbf{d}^{\mathrm{T}} \mathbf{N}^{\mathrm{T}} \mathbf{t} \mathrm{d} \Gamma=0.
\end{equation}

As the provided expression remains valid for any arbitrary virtual displacement , we can conclude that:
\begin{equation}
\int_{\Omega} \mathbf{B}^{\mathrm{T}} \boldsymbol{\sigma} \mathrm{d} \Omega-\int_{\Omega} \mathbf{N}^{\mathrm{T}}\mathbf{b} \mathrm{d} \Omega-\int_{\Gamma_t} \mathbf{N}^{\mathrm{T}} \mathbf{t} \mathrm{d} \Gamma=0.
\end{equation}

Following this, we can express the discrete governing equation as following:
\begin{equation}
\mathbf{K d}=\mathbf{P},
\end{equation}
\begin{equation}
\mathbf{K}=\int_{\Omega} \mathbf{B}^{\mathrm{T}} \mathbf{D B} \mathrm{d} \Omega,
\end{equation}
\begin{equation}
\mathbf{P}=\int_{\Omega} \mathbf{N}^{\mathrm{T}} \mathbf{b} \mathrm{d} \Omega+\int_{\Gamma_t} \mathbf{N}^{\mathrm{T}} \mathbf{t} \mathrm{d} \Gamma,
\end{equation}
where $\mathbf{K}$ is the stiffness matrices respectively; $\mathbf{P}$ is the external force vector. 

The modified strain field $\tilde{\boldsymbol{\varepsilon}}^h$ is determined through a weighted-average formulation of the conventional strain field \citep{liuSmoothedFiniteElement2007}
\begin{equation}
\tilde{\boldsymbol{\varepsilon}}^h=\int_{\Omega_C} \boldsymbol{\varepsilon}^h(\mathbf{x}) \Phi\left(\mathbf{x}-\mathbf{x}_C\right) \mathrm{d} \Omega, \label{eq:epsion_orgin}
\end{equation}
where $\Omega_C$ denotes the sub-cell, and $\Phi\left(\mathbf{x}-\mathbf{x}_C\right)$ represents the smoothing function, which is defined as:
\begin{equation}
\Phi\left(\mathbf{x}-\mathbf{x}_C\right)= \begin{cases}1 / A_C & \mathbf{x} \in \Omega_C \\ 0 & \mathbf{x} \notin \Omega_C\end{cases}, \label{eq:smoothing function}
\end{equation}
where $A_C$ is the area of the sub-cell. The modified strain field $\tilde{\boldsymbol{\varepsilon}}^h$ is computed from the smoothed discretized gradient operator $\widetilde{\mathbf{B}}$ and the nodal displacement $\mathbf{d}$. By applying the technique of integration by parts to Eq.(\ref{eq:epsion_orgin}), we obtain:
\begin{equation}
\tilde{\varepsilon}^h\left(\mathbf{x}_C\right)=\int_{\Gamma} \mathbf{u}^h(\mathbf{x}) \mathrm{n}(\mathbf{x}) \Phi\left(\mathbf{x}-\mathbf{x}_C\right)  \mathrm{d} \Gamma-\int_{\Omega} \mathbf{u}^h(\mathbf{x}) \nabla \Phi\left(\mathbf{x}-\mathbf{x}_C\right)  \mathrm{d} \Omega. \label{eq:new_epsion}
\end{equation}

Upon substituting Eq.(\ref{eq:smoothing function}) into Eq.(\ref{eq:new_epsion}), the smoothed gradient or strain field can be derived
\begin{equation}
\tilde{\boldsymbol{\varepsilon}}^h\left(\mathbf{x}_C\right)=\int_{\Gamma_C} \mathbf{u}^h(x) \mathbf{n}(x) \Phi\left(\mathbf{x}-\mathbf{x}_C\right) \mathrm{d} \Gamma=\frac{1}{A_c} \int_{\Gamma_c} \mathbf{u}^h(\mathbf{x}) \mathbf{n}(\mathbf{x}) d \Gamma,
\end{equation}
where $\Gamma_C$ represents the boundary of the sub-cell. The smoothed strain can be expressed as: 
\begin{equation}
\tilde{\boldsymbol{\varepsilon}}\left(\mathbf{x}_C\right)=\sum_I^n \widetilde{\mathbf{B}}_I\left(\mathbf{x}_C\right) \mathbf{d}_I,
\end{equation}
where $n$ denotes the number of element nodes. The stiffness matrix computed for the CSFEM is presented as follows:
\begin{equation}
\widetilde{\mathbf{K}}=\sum_{I=1}^{n e l} \widetilde{\mathbf{K}}^h=\sum_{I=1}^{n e l} \int_{\Omega^e} \widetilde{\mathbf{B}}^{\mathrm{T}} \mathbf{D} \widetilde{\mathbf{B}} \mathrm{d} \Omega. \label{eq:kmatrix}
\end{equation}

The derivation process also reveals that the strain-displacement matrix
\begin{equation}
\overset{\sim}{\mathbf{B}}=\begin{bmatrix}\overset{\sim}{\mathbf{B}}_1&\overset{\sim}{\mathbf{B}}_2&\cdots&\overset{\sim}{\mathbf{B}}_n\end{bmatrix},
\end{equation}
where the strain-displacement matrix can be written as
\begin{equation}
\overset{\sim}{\mathbf{B}}_{k}=\frac{1}{A_c}\int_{\Gamma_c^h}\begin{bmatrix}n_x(\mathbf{x})&0\\0&n_y(\mathbf{x})\\n_y(\mathbf{x})&n_x(\mathbf{x})\end{bmatrix}N_k(\mathbf{x})\mathrm{d}\Gamma. \label{eq:bk}
\end{equation}

\section{Elasto-plastic constitutive formulation}
\label{Elasto-plastic constitutive formulation}
In elasto-plastic analysis, the total incremental strain $\Delta \boldsymbol{\varepsilon}$ is conventionally decomposed into its elastic and plastic components as
\begin{equation}
\Delta \boldsymbol{\varepsilon} = 
\Delta \boldsymbol{\varepsilon}^{\mathrm{e}} + 
\Delta \boldsymbol{\varepsilon}^{\mathrm{p}},
\end{equation}
where $\Delta \boldsymbol{\varepsilon}^{\mathrm{e}}$ and $\Delta\boldsymbol{ \varepsilon}^{\mathrm{p}}$ denote the elastic and plastic strain increments, respectively. The incremental stress–strain relationship for the elastic response follows Hooke’s law, given by
\begin{equation}
\Delta \boldsymbol{\sigma} 
= \mathbf{D} \left( {\Delta \boldsymbol\varepsilon} - 
\Delta \boldsymbol{\varepsilon}^{\mathrm{p}} \right)
= \mathbf{D} \Delta\boldsymbol{ \varepsilon}^{\mathrm{e}},
\end{equation}
where $\mathbf{D}$ represents the elastic constitutive matrix.

The evolution direction of plastic strain is governed by the plastic flow rule. For associative plasticity, the plastic strain increment is determined from the yield function $F$ and the plastic multiplier $\Delta\lambda$ as
\begin{equation}
\Delta \boldsymbol{\varepsilon}^{\mathrm{p}} 
= \frac{\partial g}{\partial \boldsymbol{\sigma}} \, \Delta\lambda
= \frac{\partial F}{\partial \boldsymbol{\sigma}} \, \Delta\lambda,
\end{equation}
where $g$ is the plastic potential function, which equals the yield function $F$ in the associative case. When the stress state reaches the yield surface, the consistency condition must be satisfied to maintain stress continuity, expressed as
\begin{equation}
\Delta F = 
\frac{\partial F}{\partial \boldsymbol{\sigma}} \Delta \boldsymbol{\sigma} = 0.
\end{equation}

In this study, the Mohr–Coulomb failure criterion is employed. The yield function is defined as \citep{sherif2012matlab}
\begin{equation}
F =
\left( \sqrt{3}\cos\theta - \sin\theta \sin\phi \right) q
- 3p \sin\phi - 3c\cos\phi = 0,
\end{equation}
where $p = \tfrac{1}{3}(\sigma_1 + \sigma_2 + \sigma_3)$ is the mean stress, 
$q = \sqrt{\tfrac{3}{2} \boldsymbol{s} : \boldsymbol{s}}$ is the equivalent deviatoric stress, 
$\theta$ is the Lode angle, 
$\phi$ is the internal friction angle, and 
$c$ denote the soil cohesion, respectively.

By combining the flow rule, the consistency condition, and the hardening law, the elasto-plastic constitutive matrix can be derived as \citep{jabbarilakNonlinearScaledBoundary2019}
\begin{equation}
\mathbf{D}_{ep}=\mathbf{D}-\mathbf{D}
\left(\frac{\partial F}{\partial\boldsymbol{\sigma}}\right)
\left(\frac{\partial F}{\partial\boldsymbol{\sigma}}\right)^{\mathrm{T}}
\mathbf{D}
\left(\left(\frac{\partial F}{\partial\boldsymbol{\sigma}}\right)^{\mathrm{T}}
\mathbf{D}
\left(\frac{\partial F}{\partial\boldsymbol{\sigma}}\right)\right)^{-1}.
\label{eq:epDmatrix}
\end{equation}

The elasto-plastic constitutive matrix $\mathbf{D}_{ep}$ provides the incremental stress–strain relationship during yielding and varies spatially within the domain, reflecting the nonlinear mechanical response of the material.

In the numerical solution procedure, the equilibrium equation is expressed as
\begin{equation}
\mathbf{R} = \mathbf{F}_{\mathrm{ext}} - \mathbf{F}_{\mathrm{int}} = \mathbf{0},
\end{equation}
where $\mathbf{R}$ is the residual force vector. The internal nodal force vector is evaluated by integrating the stresses over the element domain as
\begin{equation}
\mathbf{F}_{\mathrm{int}}^{e} 
= \int_{\Omega^e} \overset{\sim}{\mathbf{B}}^{\mathrm{T}} \boldsymbol{\sigma} \, \mathrm{d}\Omega.
\end{equation}

During each iteration, the stress increment is obtained from
\begin{equation}
\Delta \boldsymbol{\sigma} = \mathbf{D}_{ep} \Delta \boldsymbol{\varepsilon},
\end{equation}
and the residual force is updated accordingly. The nonlinear system is solved using a Newton–Raphson scheme:
\begin{equation}
\mathbf{K}_{t}^{(i)} \Delta \mathbf{u}^{(i)} = \mathbf{R}^{(i)},
\end{equation}
where $\mathbf{K}_{t}^{(i)}$ is the current tangent stiffness matrix. Iterations proceed until the normalized residual force satisfies
\begin{equation}
\frac{\|\mathbf{R}^{(i)}\|}{\|\mathbf{F}_{\mathrm{ext}}\|} \leq \varepsilon_{\mathrm{tol}}.
\end{equation}

This residual force evaluation ensures equilibrium at each load step and enables accurate simulation of nonlinear elasto-plastic deformation.

\section{Construction of CSFEM elements}
\label{Construction of CSFEM elements}
\begin{figure}[H]
  \centering
  \includegraphics[width=1.0\textwidth]{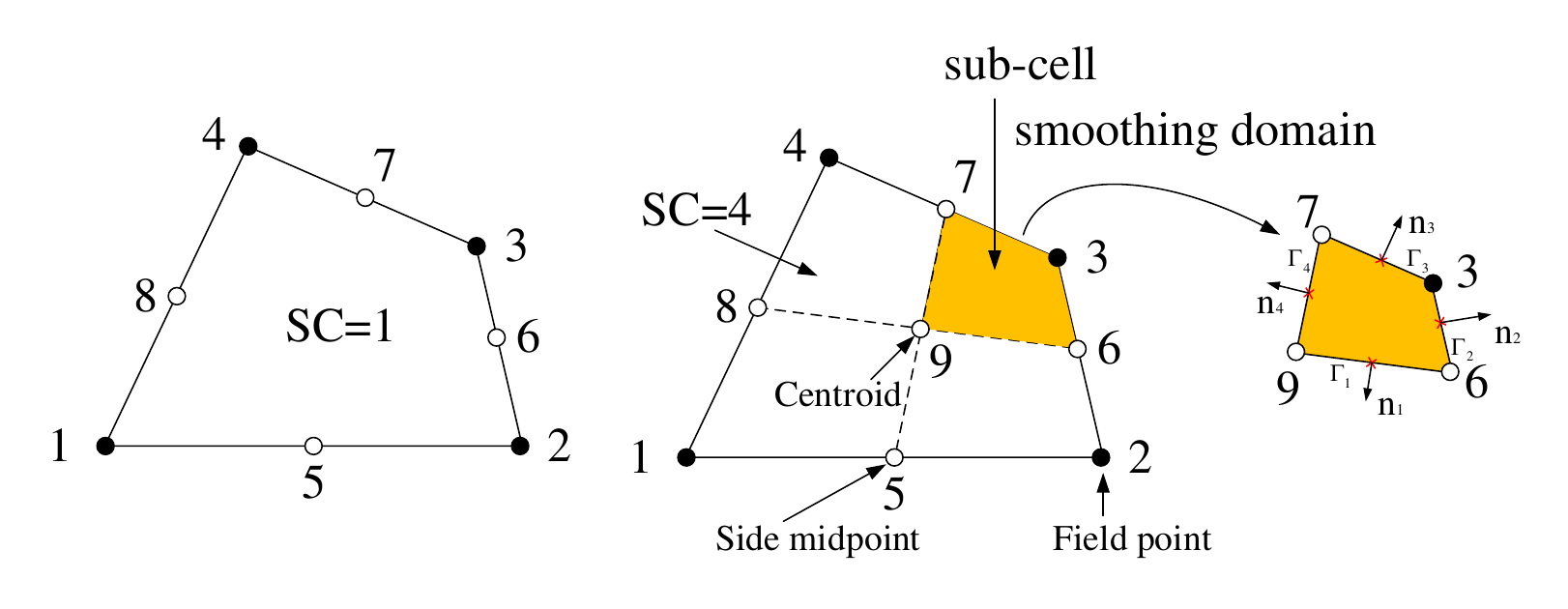}
  \caption{Construction of a cell-based smoothing domain.}
  \label{fig:csfem_element}
\end{figure}

In the CSFEM, each element is subdivided into several smoothing domains to construct the stiffness matrix. As shown in Fig.~\ref{fig:csfem_element}, a quadrilateral element is divided into four subcells, within which the strain field is smoothed and the domain integration is transformed into boundary integrals. Since only the shape function values along the boundaries are required, simple interpolation schemes are used at the Gauss points located on the element edges. Representative shape function values at typical points within the element are listed in Tab.~\ref{tab:shape_function_values}, which are used to evaluate the smoothed strain field and the corresponding stiffness matrix.

\begin{table}[H]
  \centering
  \caption{Shape function value at different sites within an element (see Fig.~\ref{fig:csfem_element}).}
  \label{tab:shape_function_values}
  \begin{tabular}{cccccc}
    \toprule
    Site & Node 1 & Node 2 & Node 3 & Node 4 & Description \\
    \midrule
    1 & 1.0  & 0    & 0    & 0    & Field node \\
    2 & 0    & 1.0  & 0    & 0    & Field node \\
    3 & 0    & 0    & 1.0  & 0    & Field node \\
    4 & 0    & 0    & 0    & 1.0  & Field node \\
    5 & 0.5  & 0.5  & 0    & 0    & Side midpoint \\
    6 & 0    & 0.5  & 0.5  & 0    & Side midpoint \\
    7 & 0    & 0    & 0.5  & 0.5  & Side midpoint \\
    8 & 0.5  & 0    & 0    & 0.5  & Side midpoint \\
    9 & 0.25 & 0.25 & 0.25 & 0.25 & Centroid\\
    \bottomrule
  \end{tabular}
  \label{tab:shape_function_value}
\end{table}

The use of smoothed strain fields in CSFEM leads to a softer stiffness matrix compared with the conventional FEM, which helps alleviate numerical issues such as mesh distortion sensitivity and volumetric locking. Furthermore, the method retains simplicity in implementation since only the original finite element shape functions are required at the boundary integration points, without computing their spatial derivatives.

\section{Implementation}
\label{Implementation}
\subsection{Implementation of the elasto-plastic CSFEM}
\begin{figure}[H]
  \centering
  \includegraphics[width=1.0\textwidth]{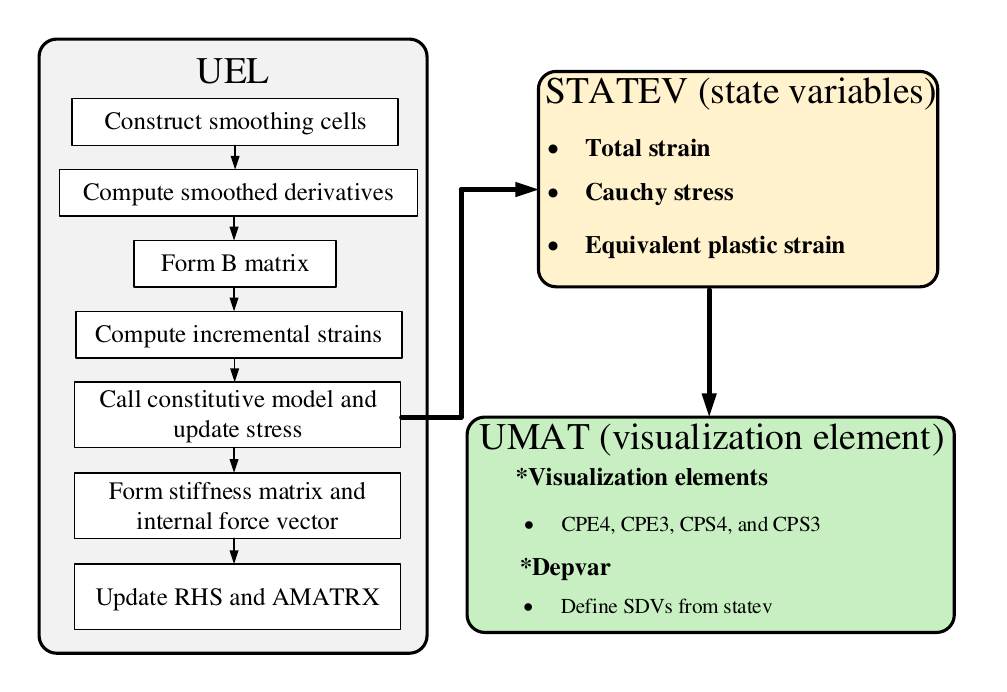}
  \caption{Framework of implementing the elasto-plastic CSFEM in ABAQUS using UEL
and UMAT subroutines.}
  \label{fig:UEL2UMAT}
\end{figure}

Fig.~\ref{fig:UEL2UMAT} and Algorithm~\ref{alg1} illustrate the computational workflow of the elastic--plastic analysis implemented in ABAQUS via the UEL interface. The main responsibility of the UEL subroutine is to evaluate and update the element contributions to the internal force vector (RHS) and stiffness matrix (AMATRX), based on the information supplied by the ABAQUS/Standard solver~\citep{abaqus2010user}. Within this framework, the stiffness matrix and the residual force vector are formulated as
\begin{equation}\mathrm{AMATRX}=\tilde{\mathrm{K}},\label{eq:amatrx}\end{equation}
\begin{equation}\mathrm{RHS}=-\tilde{\mathrm{K}}\mathrm{U}_{t+\Delta t},\label{eq:rhs}\end{equation}
where $\tilde{\mathrm{K}}$ denotes the element stiffness matrix and $\mathrm{U}_{t+\Delta t}$ represents the displacement vector at the end of the current increment. The UEL routine iteratively updates both AMATRX and RHS according to Eqs.~(\ref{eq:amatrx}) and~(\ref{eq:rhs}) throughout the solution process.

\begin{algorithm}[H]
\caption{Solving the elasto-plastic CSFEM UEL}
\begin{algorithmic}[1]
\Require Node and element information, material properties, nodal displacement $\mathbf{u}_t$
\Ensure Nodal displacement $\mathbf{u}_{t+1}$
\State \textbf{Iteration:} increment steps $k = 1$
\While{ABAQUS not converged}
    \State solve nodal displacement $\mathbf{u}^{k}_{t+1}$
    \For{$i = 1$ to $n_c$}
        \State solve shape functions $N_k(x)$ in Tab.\ref{tab:shape_function_value}
        \State solve strain–displacement matrix $\tilde{\mathbf{B}}_c$ in Eq.(\ref{eq:bk})
        \State solve elasto-plastic constitutive matrix $\mathbf{D}_{ep}$ in Eq.(\ref{eq:epDmatrix})
        \State solve stiffness matrix of subcells $\tilde{\mathbf{K}}_c$ in Eq.(\ref{eq:kmatrix})
    \EndFor
    \State obtain element matrices $\tilde{\mathbf{K}}$ 
    \State update AMATRX and RHS according to Eqs.~(\ref{eq:amatrx}) and (\ref{eq:rhs})
    \State update $k = k + 1$
\EndWhile
\State solve nodal field variables $\mathbf{u}_{t+1} = \mathbf{u}^{k}_{t+1}$
\end{algorithmic}
\label{alg1}
\end{algorithm}

\subsection{UEL–UMAT data transfer mechanism for stress and strain visualization}
\begin{figure}[H]
  \centering
  \includegraphics[width=1.0\textwidth]{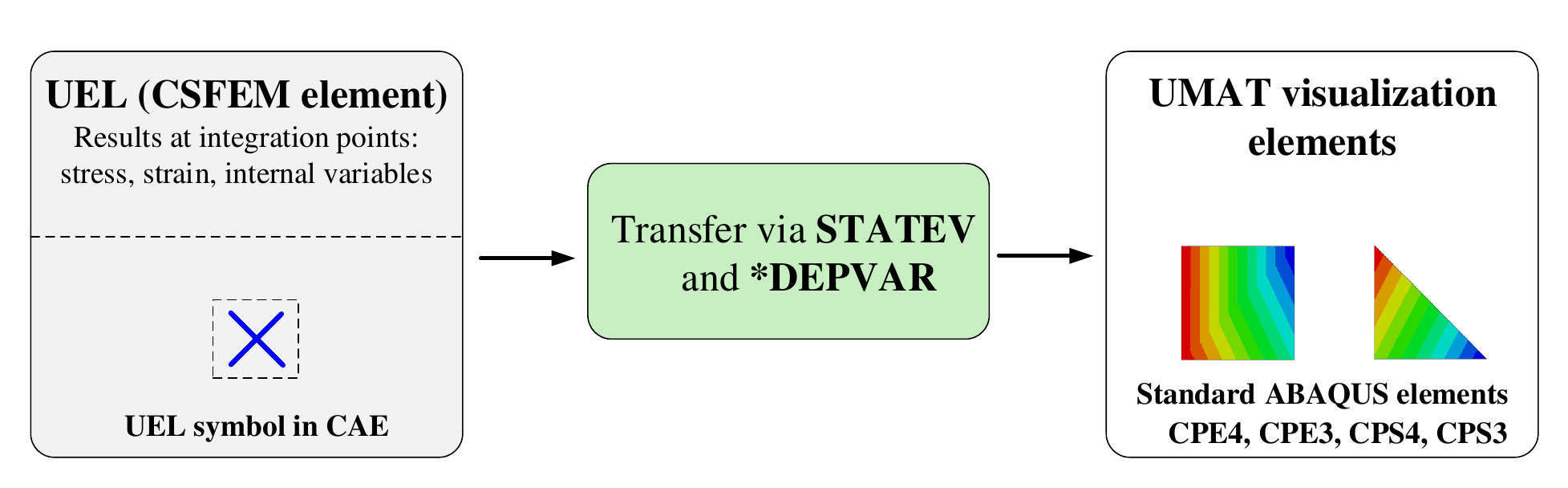}
  \caption{Data-linking scheme for transferring stress, strain, and other state variables from the UEL to UMAT-based visualization elements in ABAQUS.}
  \label{fig:visual_element}
\end{figure}

As shown in Fig.~\ref{fig:visual_element}, the UEL in ABAQUS is represented by an “×” symbol, and its stress–strain results cannot be visualized directly. To address this limitation, a data-linking strategy was developed between the UEL and the user material subroutine (UMAT), as schematically illustrated in Fig.~\ref{fig:visual_element}. Through this interface, the stress, strain, and other state variables (STATEV) computed within the UEL were transferred to the UMAT framework, allowing the use of standard ABAQUS elements for post-processing and visualization.

For visualization in ABAQUS/CAE, these quantities were exported via state variables declared in the input file using the *DEPVAR option and were stored at the integration points during the analysis. Consequently, as shown in Fig.~\ref{fig:visual_element}, the numerical results obtained from the UEL could be accessed and visualized within the standard post-processing environment, providing clear field distributions comparable to those produced by conventional built-in elements.

To facilitate the transfer of stress, strain, and internal variables from the UEL to the UMAT framework, the state variables were explicitly declared in the input file through the *DEPVAR option. The corresponding ABAQUS syntax is listed below:

\begin{lstlisting}[caption={ABAQUS input for defining state variables via *DEPVAR},label={lst:depvar}]
*Material, name=Material-1
*User Material, constants=0
*Depvar
11,
1, EPSI, EPSI
2, TEPS, TEPS
3, PE11, PE11
4, PE22, PE22
5, PE33, PE33
6, PE12, PE12
7, S11, S11
8, S22, S22
9, S33, S33
10, S12, S12
11, PEMAG, PEMAG
\end{lstlisting}

The above definitions specify 11 state variables that store the quantities computed within the UEL at each 
integration point. For clarity, the physical meaning of these variables is summarized in 
Tab.~\ref{tab:depvar}.

\begin{table}[H]
\centering
\caption{Definition of state variables used for UEL--UMAT data transfer.}
\begin{tabular}{ccc}
\toprule
Index & Variable name & Description \\
\midrule
1  & EPSI   & Equivalent plastic strain increment \( \Delta \bar{\varepsilon}^{p} \)\\
2  & TEPS   & Total strain measure \\
3  & PE11   & Plastic strain component \( \varepsilon^{p}_{11} \) \\
4  & PE22   & Plastic strain component \( \varepsilon^{p}_{22} \) \\
5  & PE33   & Plastic strain component \( \varepsilon^{p}_{33} \) \\
6  & PE12   & Plastic strain component \( \varepsilon^{p}_{12} \) \\
7  & S11    & Stress component \( \sigma_{11} \) \\
8 & S22    & Stress component \( \sigma_{22} \) \\
9 & S33    & Stress component \( \sigma_{33} \) \\
10 & S12    & Stress component \( \sigma_{12} \) \\
11 & PEMAG  & Total equivalent plastic strain \( \bar{\varepsilon}^{p} \) \\
\bottomrule
\end{tabular}
\label{tab:depvar}
\end{table}

\section{Numerical examples}
\label{Numerical examples}
This section introduces several benchmark examples designed to verify the accuracy and convergence characteristics of the proposed elasto-plastic analysis framework. The numerical responses obtained with the CSFEM formulation are systematically compared with reference solutions generated using standard finite element analyses in ABAQUS. All computations were performed on a workstation equipped with an Intel Core i7-4710MQ processor (2.50~GHz) and 4~GB of RAM. The quality of the numerical results is quantified using the following relative error measure:
\begin{equation}
    \mathbf{e}_{L_2}=\frac{\left\|\mathbf{U}_{num}-\mathbf{U}_{ref}\right\|_{L_2}}{\parallel\mathbf{U}_{ref}\parallel_{L_2}},
\end{equation}
where $\mathbf{U}_{num}$ represents the numerical solution, and $\mathbf{U}_{ref}$ denotes the reference solution.

\subsection{Pressurized thick-cylinder modeled using a quarter-annulus model}
\begin{figure}[H]
  \centering
  \includegraphics[width=0.5\textwidth]{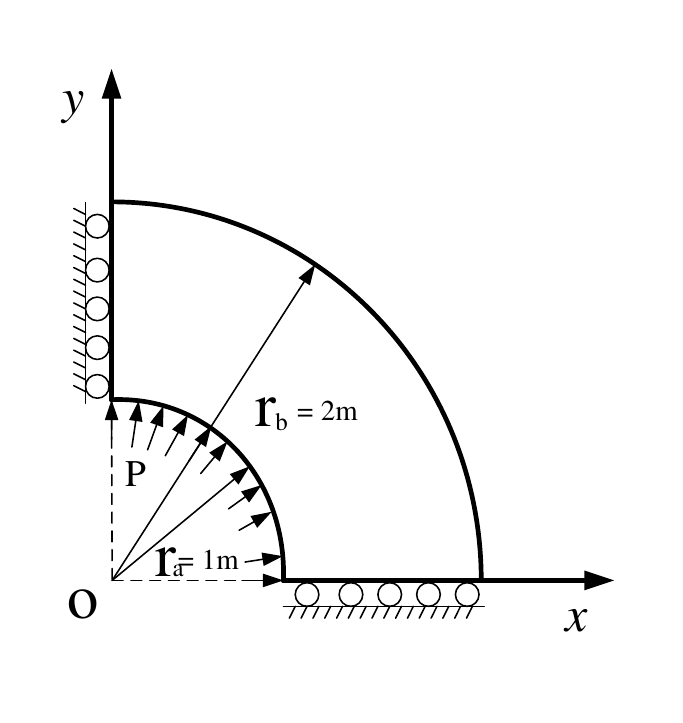}
\caption{Quarter model of a thick-walled cylinder subjected to internal pressure at the inner surface.}
  \label{fig:ex01_geo}
\end{figure}

\begin{figure}[H]
  \centering
  \includegraphics[width=0.6\textwidth]{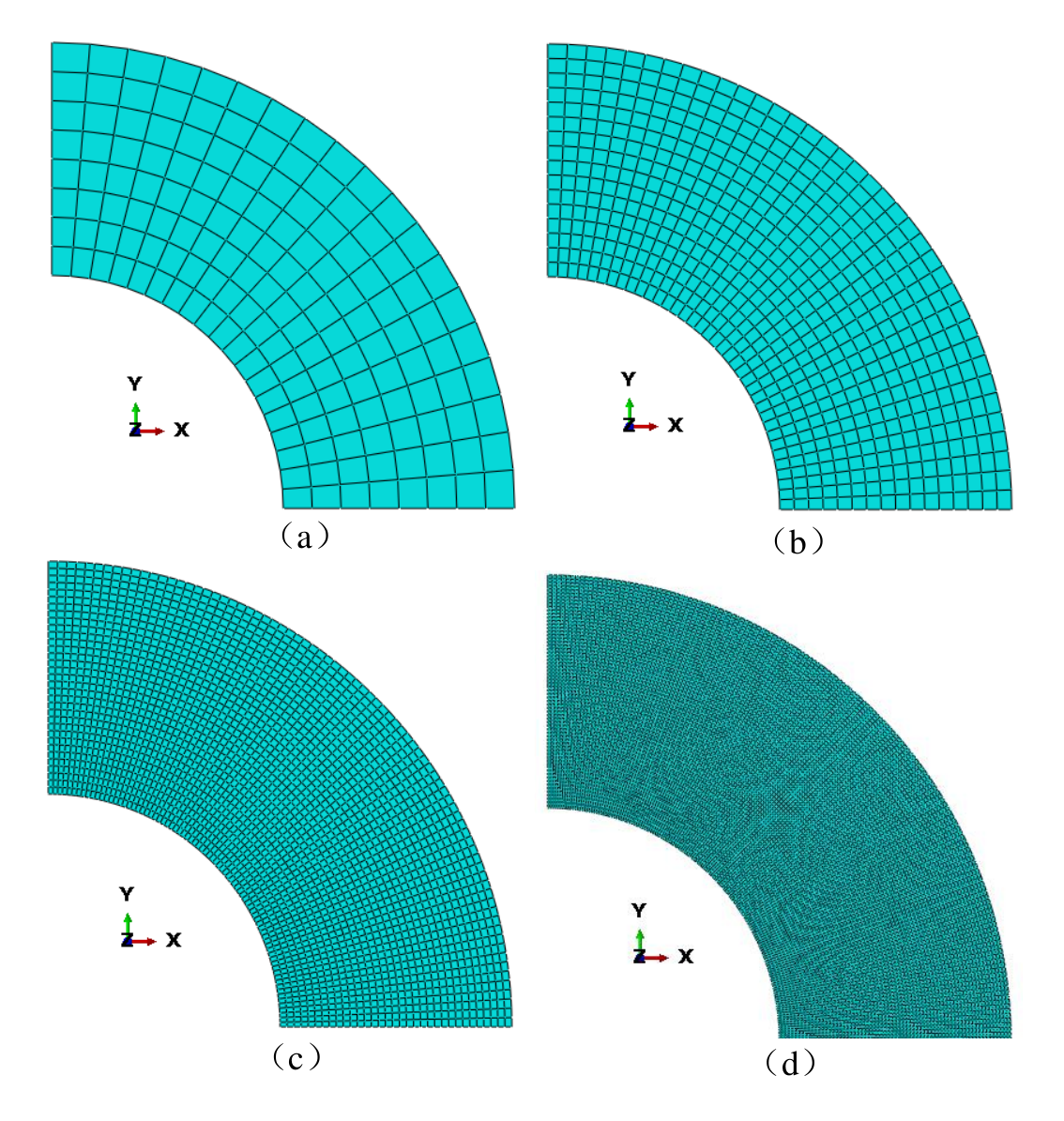}
  \caption{Four levels of mesh refinement; (a) 0.125 m mesh size; (b) 0.0625 m mesh size; (c) 0.03 m mesh size; (d) 0.015 m mesh size.}
  \label{fig:ex01_mesh}
\end{figure}

In this example, we considered a quarter-thickness cylindrical model subjected to internal pressure to verify the proposed method. The displacements in the polar coordinate system can be derived analytically~\citep{Timoshenko1970}. 
\begin{equation}
    {U}_{\mathrm{rad}}=\frac{R_a^2Pr}{E\left(R_b^2-R_a^2\right)}\left(1-\nu+\left(\frac{R_b}{r}\right)^2(1+\nu)\right),
\end{equation}

\begin{equation}
    u_x=u_{\mathrm{rad}}\cos\theta,
\end{equation}

\begin{equation}
    u_y=u_{\mathrm{rad}}\sin\theta.
\end{equation}
where $R_a$ and $R_b$ denote the inner and outer radii of the cylinder, respectively, with $R_a = 1\,\text{m}$ and $R_b = 2\,\text{m}$. A uniform pressure of $1000\,\text{Pa}$ was applied along the inner circular boundary, while the material parameters were defined by Young's modulus $E=10000$ Pa and Poisson's ratio $\nu=0.25$. The geometry of the cylinder and the associated boundary conditions are illustrated in Fig.~\ref{fig:ex01_geo}. The meshes are successively refined to element sizes of 0.125 m, 0.0625 m, 0.03 m, and 0.015 m, as shown in Fig.~\ref{fig:ex01_mesh}.

\begin{figure}[H]
    \centering
    \begin{subfigure}{0.48\textwidth}
        \centering
        \includegraphics[width=\linewidth]{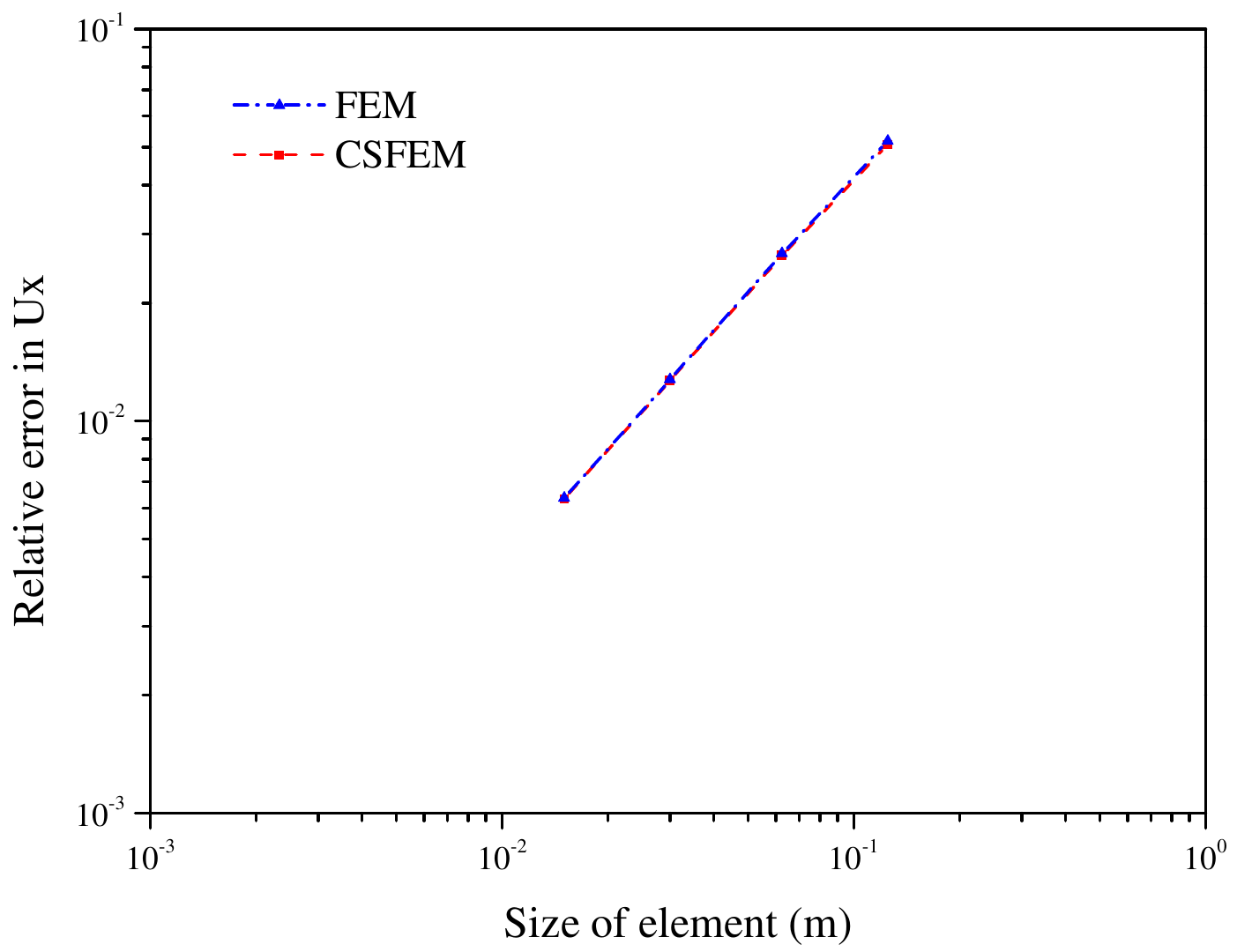}
        \caption{Relative error in $\mathrm{U_x}$}
        \label{fig:ex01_cong_Ux}
    \end{subfigure}
    \hfill
    \begin{subfigure}{0.48\textwidth}
        \centering
        \includegraphics[width=\linewidth]{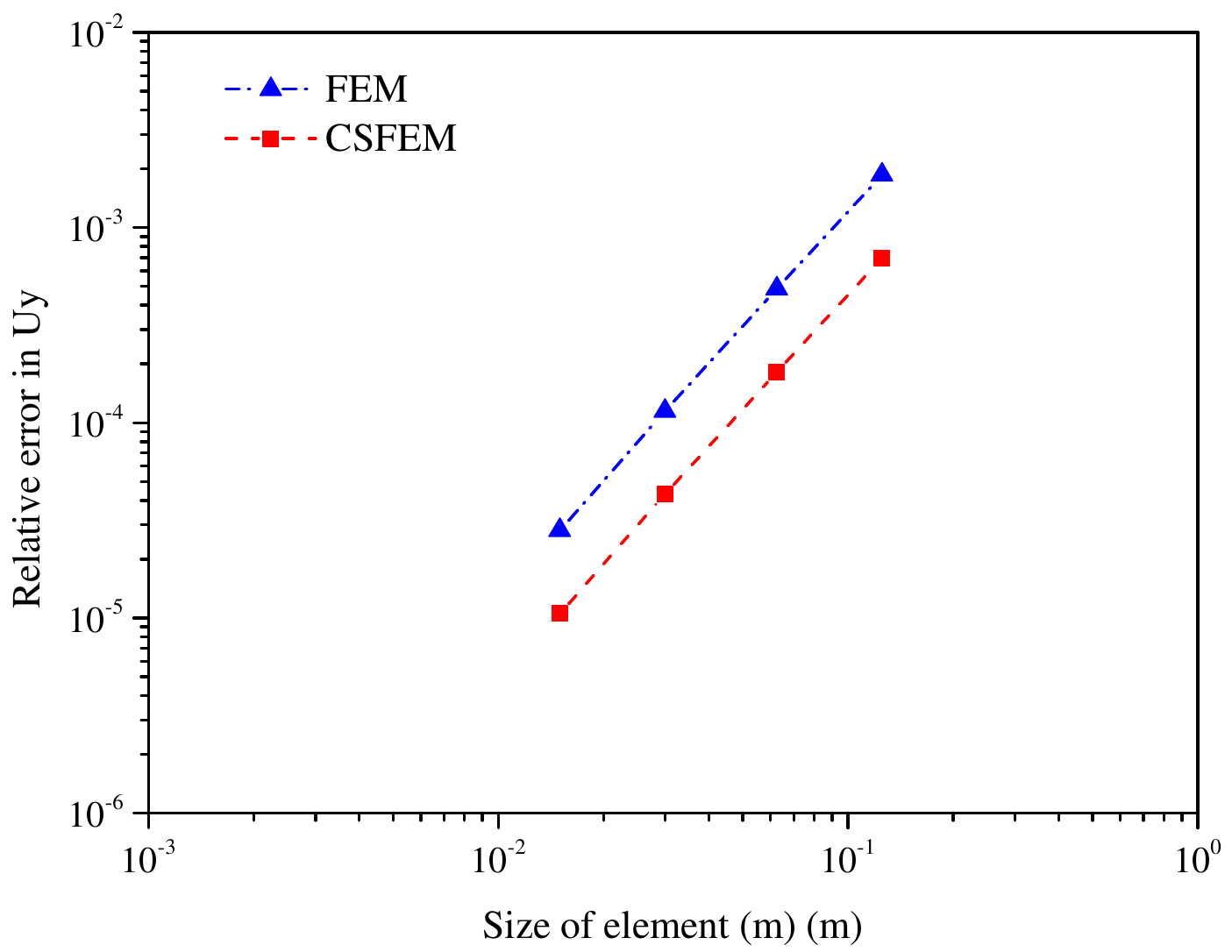}
        \caption{Relative error in $\mathrm{U_y}$  }
        \label{fig:ex01_cong_Uy}
    \end{subfigure}
    \caption{Convergence of displacement errors in the (a) horizontal direction; (b) vertical direction.}
    \label{fig:ex01_cong}
\end{figure}

\begin{figure}[H]
  \centering
  \includegraphics[width=0.7\textwidth]{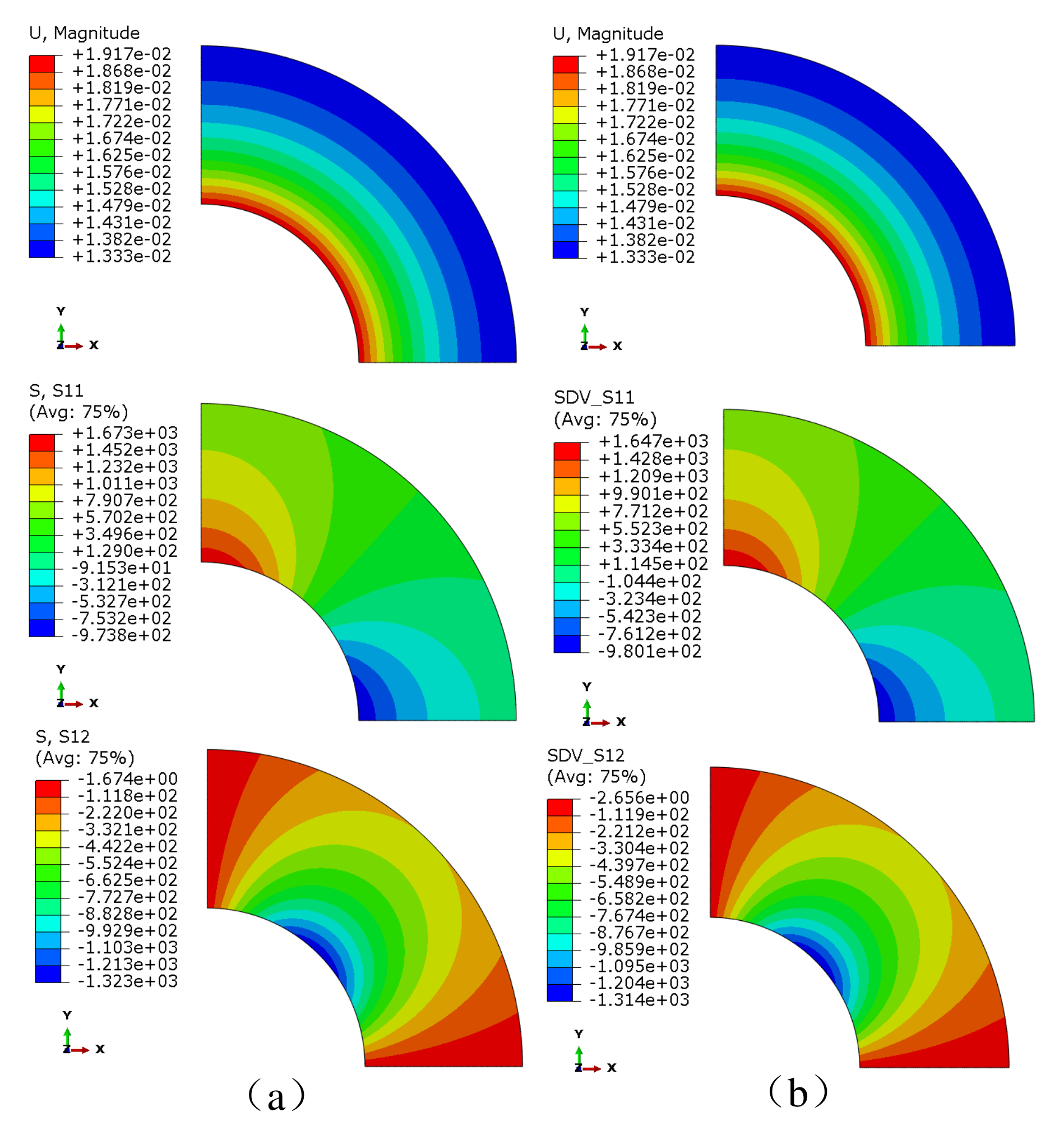}
  \caption{Contour plots of the numerical results for the thick cylinder problem: (a) results of the FEM; (b) results of the CSFEM.}
  \label{fig:ex01_contour}
\end{figure}

Fig.~\ref{fig:ex01_cong} demonstrates that both the CSFEM and FEM exhibit satisfactory convergence behavior. As summarized in Tab.~\ref{ex01:tab2}, the solution accuracy consistently improves with mesh refinement. In particular, for the finest mesh with a characteristic size of 0.003,m, the relative errors in the vertical displacement are $4.290\times10^{-5}$ for the CSFEM and $1.144\times10^{-4}$ for the FEM. Comparable trends are observed across all mesh levels, where the CSFEM yields smaller errors than the FEM. This indicates that the CSFEM offers enhanced accuracy relative to the conventional FEM. Furthermore, the displacement contour plots shown in Fig.~\ref{fig:ex01_contour} reveal an excellent agreement between the two methods, confirming the consistency of their numerical predictions.

\begin{table}[H]
    \centering
    \caption{Relative errors in the displacement with mesh refining.}
    \label{tab:relerr_disp}
    \begin{tabular}{cccccc}
        \toprule
        \multirow{3}{*}{Size of element (m)} 
            & \multicolumn{4}{c}{Relative errors} \\
        \cmidrule(lr){2-5}
            & \multicolumn{2}{c}{Ux} 
            & \multicolumn{2}{c}{Uy} \\
        \cmidrule(lr){2-3}\cmidrule(lr){4-5}
            & CSFEM & FEM & CSFEM & FEM \\
        \midrule
        0.125 & $5.07\times10^{-2}$ & $5.186\times10^{-2}$ & $6.989\times10^{-4}$ & $1.860\times10^{-3}$ \\
        0.0625 & $2.650\times10^{-2}$ & $2.680\times10^{-2}$ & $1.817\times10^{-4}$ & $4.846\times10^{-4}$ \\
        0.03  & $1.270\times10^{-2}$ & $1.277\times10^{-2}$ & $4.290\times10^{-5}$ & $1.144\times10^{-4}$ \\
        0.015 & $6.340\times10^{-3}$ & $6.360\times10^{-3}$ & $1.051\times10^{-5}$ & $2.804\times10^{-5}$ \\
        \bottomrule
    \end{tabular}
    \label{ex01:tab2}
\end{table}

\subsection{Biaxial test}
\begin{figure}[H]
  \centering
  \includegraphics[width=0.9\textwidth]{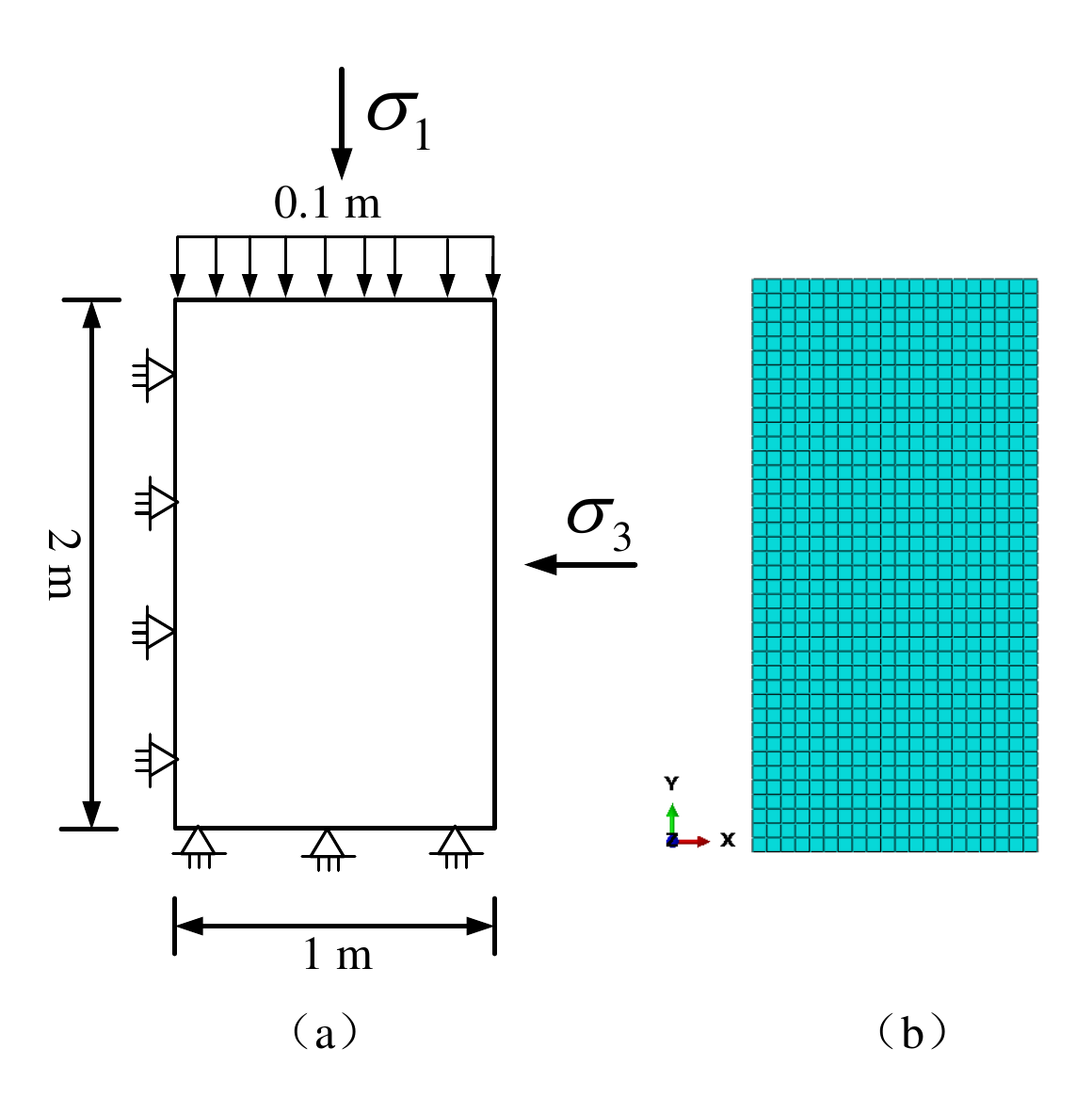}
  \caption{Schematic diagram for the biaxial test; (a) geometry and boundary conditions; (b) mesh model.}
  \label{fig:ex05_geo_mesh}
\end{figure}

In this example, we consider a typical biaxial test setup used to investigate the mechanical behavior of soils, as shown in Fig.~\ref{fig:ex05_geo_mesh} (a). The soil specimen is initially subjected to prescribed confining stresses of 100 kPa in both the vertical and horizontal directions. After the initial stress state is established, a vertical displacement is incrementally imposed while the lateral stress is kept constant until failure occurs. The soil is modeled with an elastic modulus $E = 10$ MPa, Poisson’s ratio $\nu = 0.3$, cohesion $c = 10$ kPa, and an internal friction angle of $30^{\circ}$.

For the present biaxial test, adopting the Mohr-Coulomb failure criterion with cohesion $c$ and friction angle $\phi$, 
and taking the major and minor principal stresses as $\sigma_1$ and $\sigma_3$, respectively, 
the ultimate major principal stress is given by \citep{FeiPeng2017AbaqusGeo}
\begin{equation}
 \sigma_{1,\mathrm{ult}} \;=\; \sigma_3\;\frac{1+\sin\phi}{1-\sin\phi}
 \;+\;\frac{2c\cos\phi}{1-\sin\phi}.
 \label{eq:sigma_ult}
\end{equation}

\begin{figure}[H]
  \centering
  \includegraphics[width=0.9\textwidth]{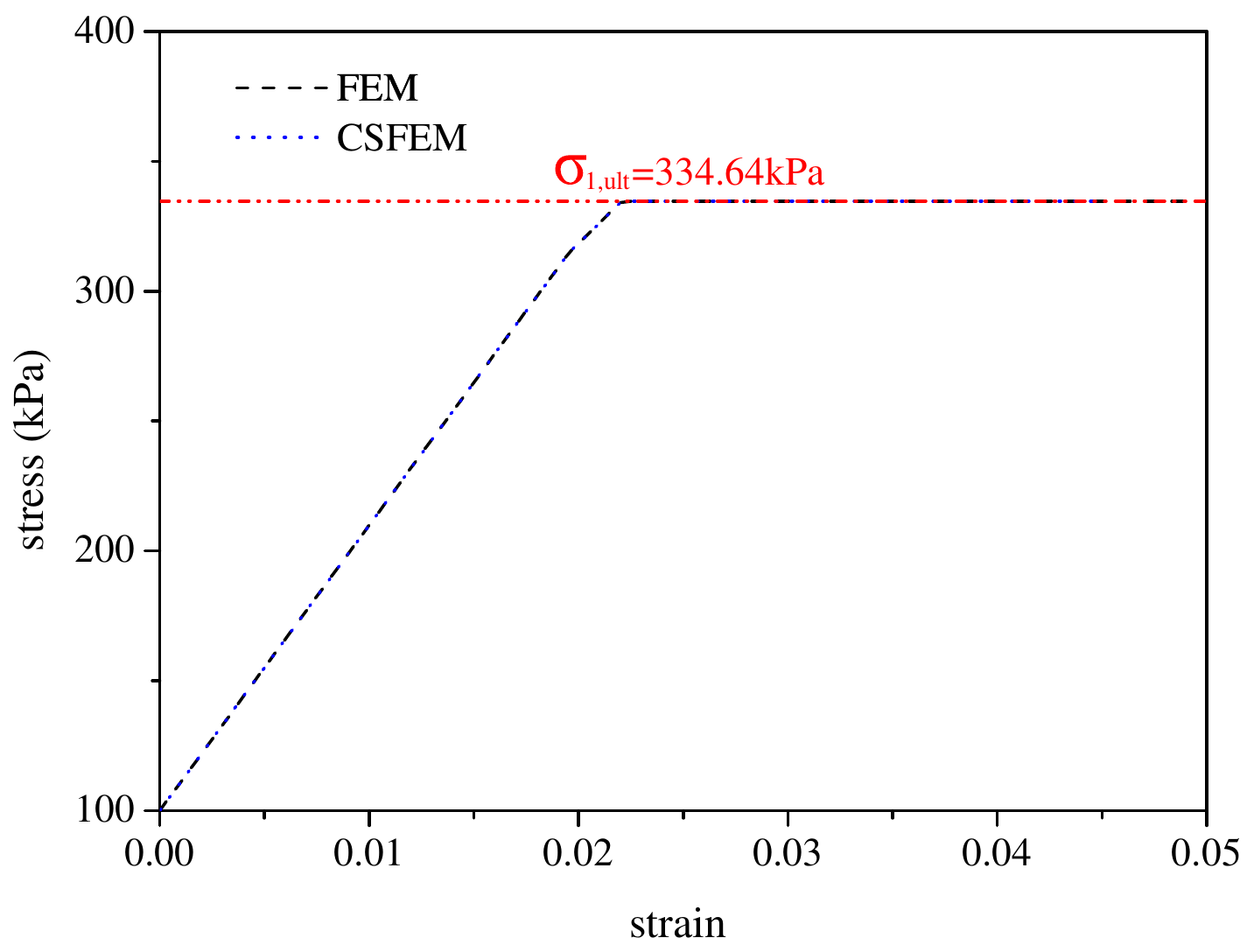}
  \caption{The relationship between vertical stress and vertical strain.}
  \label{fig:ex05_stress_strain}
\end{figure}

\begin{figure}[H]
  \centering
  \includegraphics[width=0.9\textwidth]{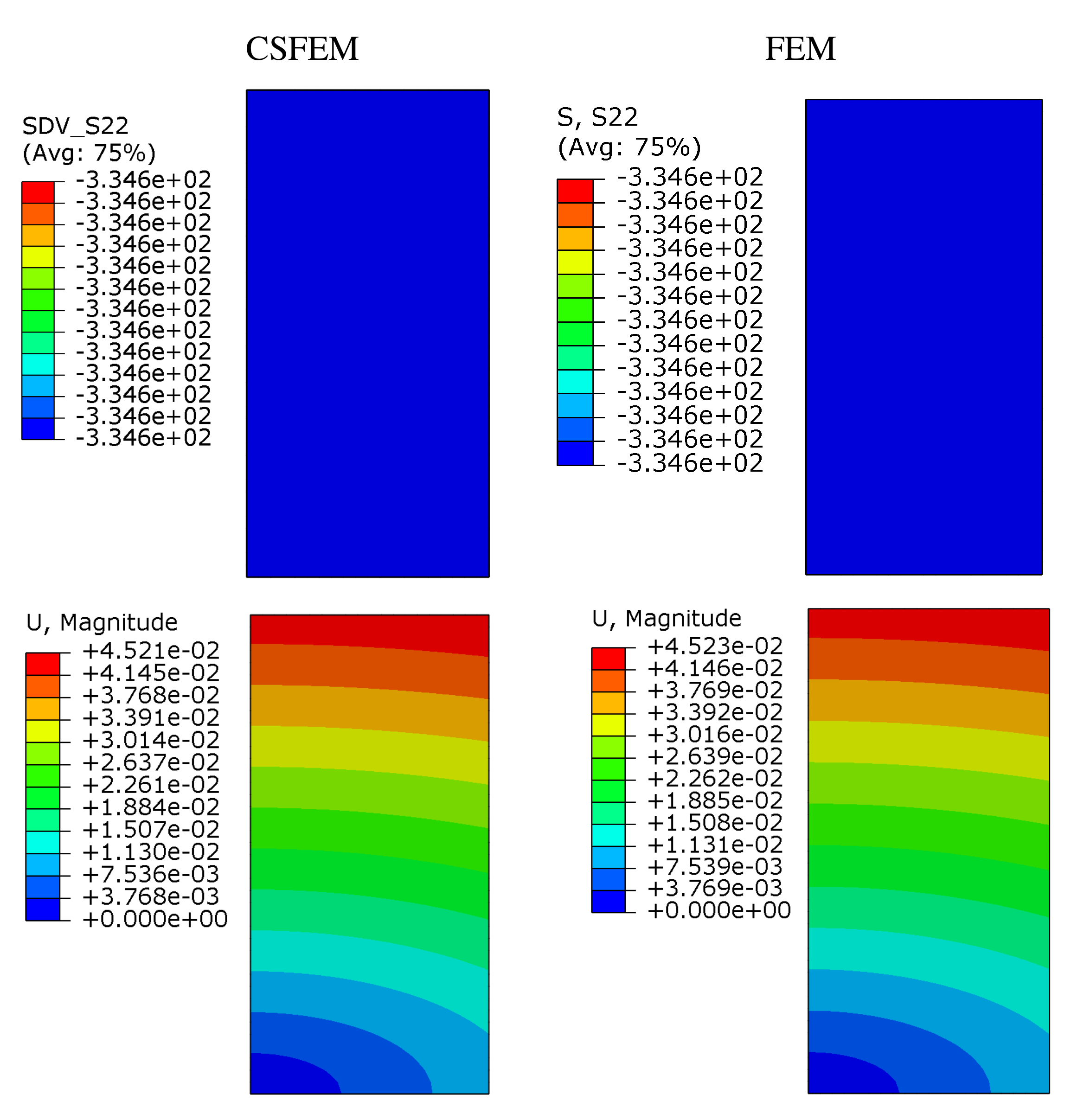}
  \caption{The vertical stress and displacement contours at the material yielding.}
  \label{fig:ex05_contour}
\end{figure}

According to Eq.~(\ref{eq:sigma_ult}), the ultimate major principal stress in this example is 334.64 kPa. Fig.~\ref{fig:ex05_stress_strain} shows the relationship between vertical stress and vertical strain. As illustrated in Fig.~\ref{fig:ex05_stress_strain}, during the initial loading stage, the material remains in the elastic region, and the vertical stress increases with the increase in vertical strain until yielding occurs. Fig.~\ref{fig:ex05_contour} illustrates the vertical stress and displacement contours at the material yielding. It can be observed from Fig.~\ref{fig:ex05_contour} that the final vertical stress obtained by CSFEM and FEM is 334.6 kPa, which agrees well with the result of Eq.~(\ref{eq:sigma_ult}). Moreover, the displacement distributions of CSFEM and FEM are also essentially consistent.

\subsection{Bearing capacity problem}
\begin{figure}[H]
  \centering
  \includegraphics[width=0.9\textwidth]{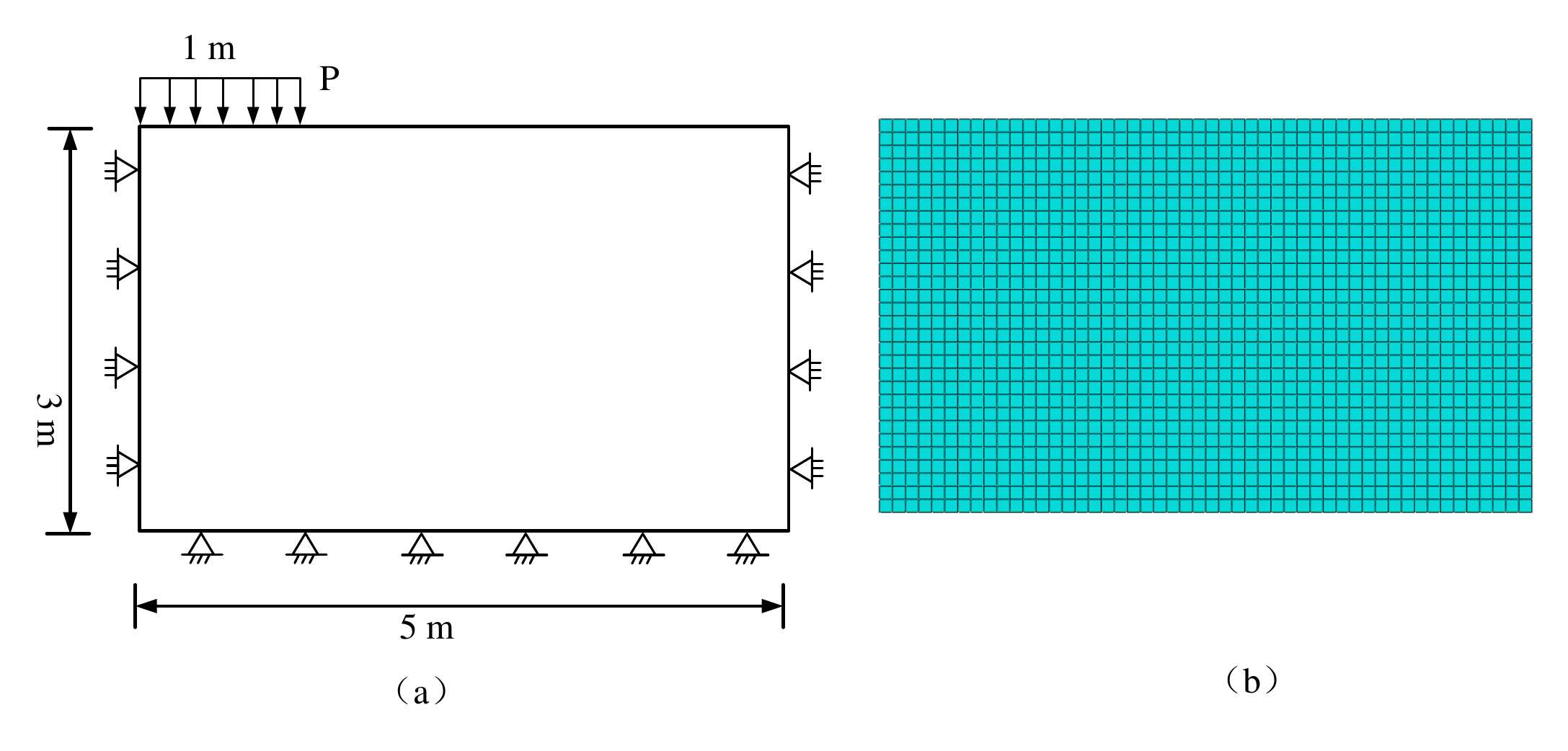}
  \caption{Schematic diagram for the bearing capacity problem; (a) geometry and boundary conditions; (b) mesh model.}
  \label{fig:ex02_geo_mesh}
\end{figure}

In this example, we consider a flexible strip footing resting on the surface of a half space subjected to a uniform vertical load, as shown in Fig.~\ref{fig:ex02_geo_mesh}(a). Due to symmetry, only half of the footing is modeled. Roller boundary conditions are imposed on the side edges, while the bottom boundary is fully fixed. The analysis is carried out under plane strain conditions using the mesh shown in Fig.~\ref{fig:ex02_geo_mesh}(b), with 80 load increments of 8000~Pa. The corresponding soil properties are given in Tab.~\ref{tab:ex02:soil_parameters}. For a soil with zero unit weight ($\gamma = 0$), the ultimate bearing capacity of strip foundations can be expressed as \citep{sherif2012matlab}:

\begin{equation}
q_u=\begin{cases}
c\cdot N_c & \phi \neq 0 \\
(\pi+2)\cdot c & \phi = 0
\end{cases} \quad,
\label{eq:q}
\end{equation}
where $N_c$ is a bearing capacity factor can be given by,
\begin{equation}
N_c=(N_q-1)\cot\phi,
\end{equation}

\begin{equation}
N_q=e^{\pi\tan\phi}\left(\frac{1+\sin\phi}{1-\sin\phi}\right).
\end{equation}

\begin{table}[H]
  \centering
  \caption{Soil parameters of a flexible strip footing.}
  \label{tab:ex02:soil_parameters}
  \begin{tabular}{ccccc}
    \toprule
    $E$ (kPa) & $\nu$ & $\phi$ & $c$ (kPa) & $\psi$ \\ 
    \midrule
    10000 & 0.3 & $5^{\circ}$  & 1 &  $5^{\circ}$ \\ 
    \bottomrule
    
  \vspace{2mm}
  {\footnotesize Note: $\psi$ denotes the dilation angle.}
  \end{tabular}
\end{table}

\begin{figure}[H]
  \centering
  \includegraphics[width=0.9\textwidth]{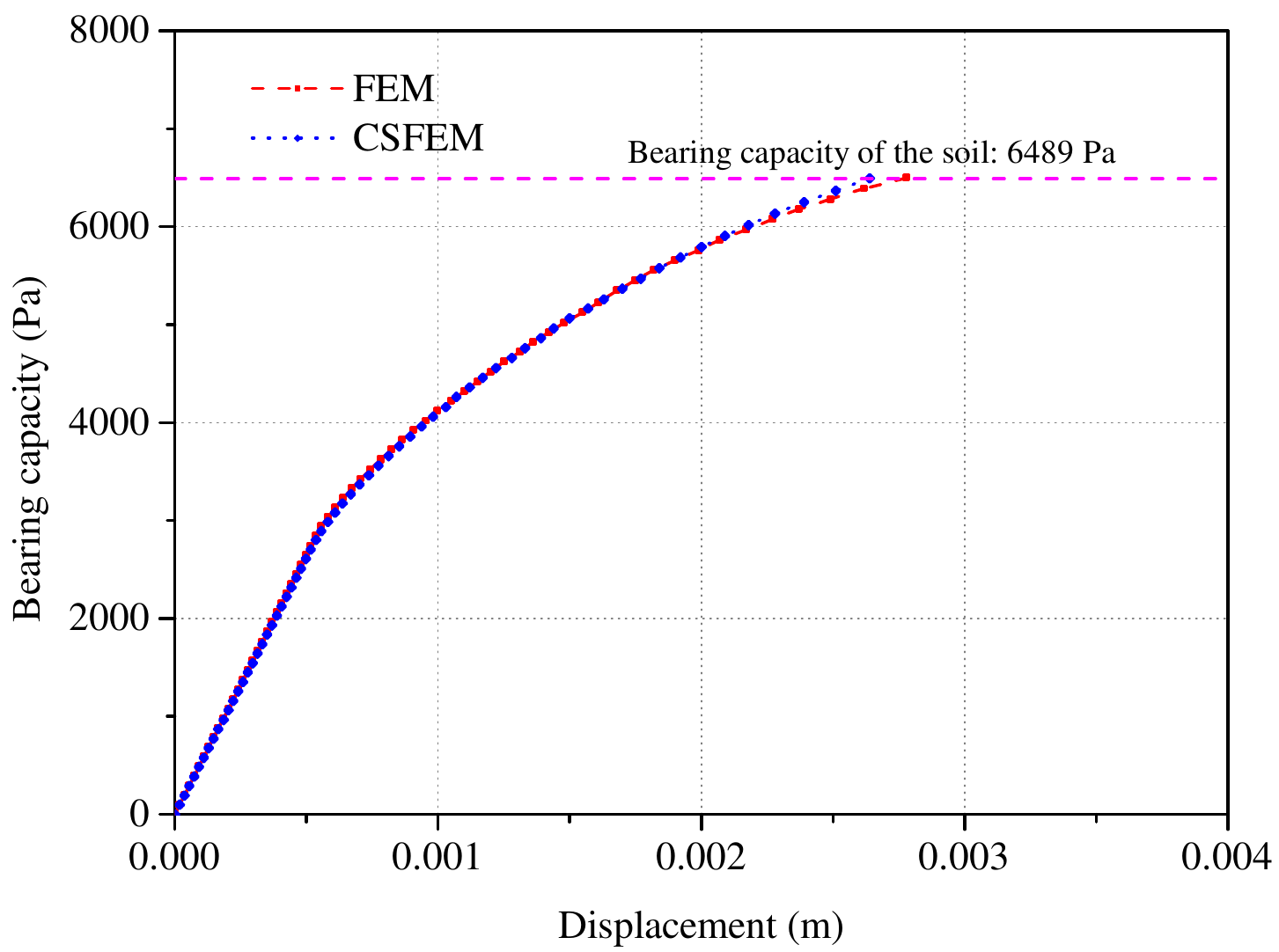}
  \caption{Bearing capacity–displacement diagram of a flexible footing using CSFEM and FEM.}
  \label{fig:ex02_disp_s11}
\end{figure}

\begin{figure}[H]
  \centering
  \includegraphics[width=1.0\textwidth]{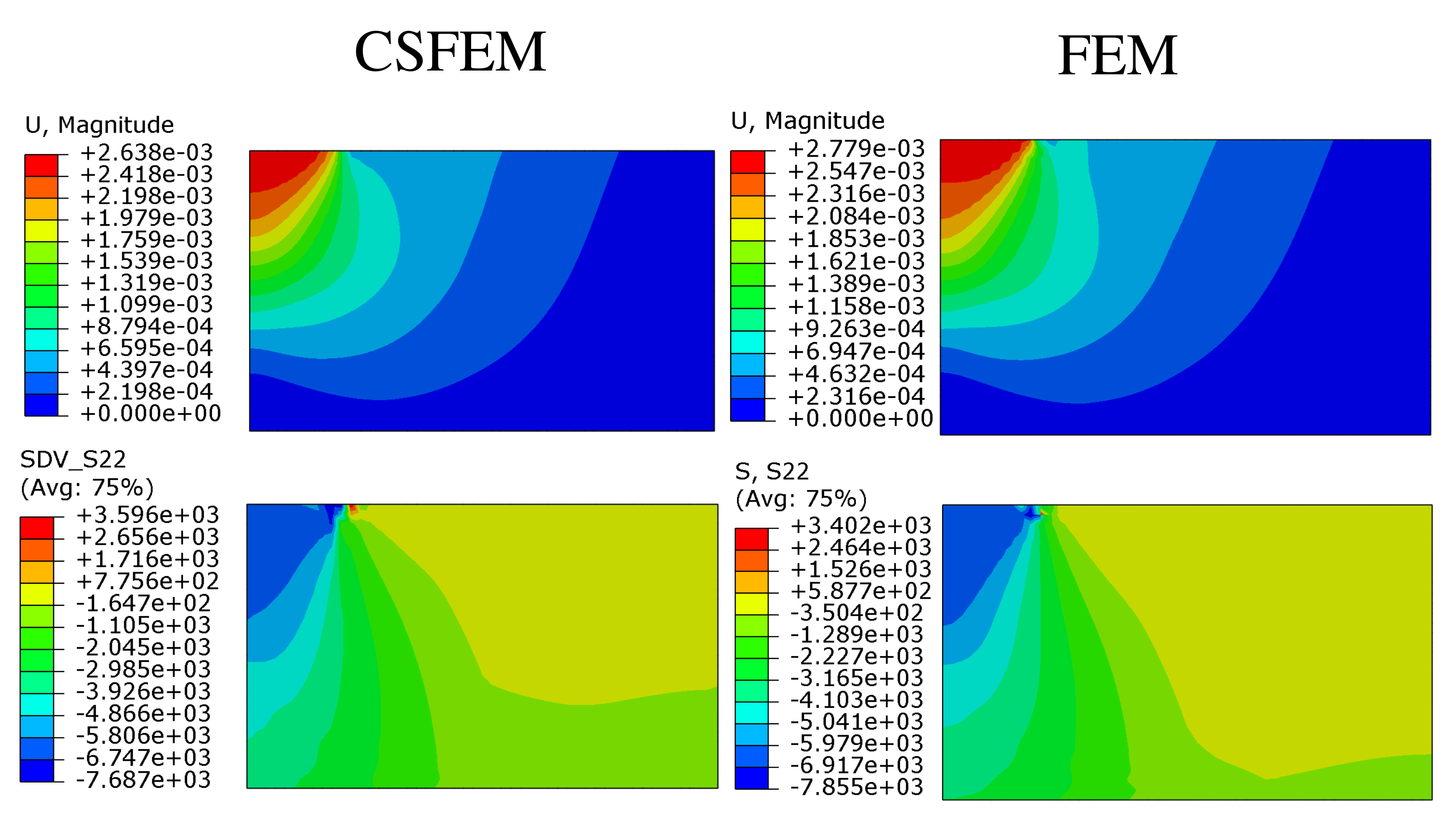}
  \caption{Comparison of vertical displacement and stress contours between CSFEM and FEM at the ultimate bearing capacity.}
  \label{fig:ex02_contour}
\end{figure}

Fig.~\ref{fig:ex02_disp_s11} presents the displacement–bearing capacity curve, indicating that the CSFEM and FEM results are in good agreement. According to Eq.(~\ref{eq:q}), the bearing capacity of the foundation in this example is 6489 Pa. Tab.~\ref{tab:ex02:Relative errors} shows that the bearing capacities obtained from CSFEM and FEM are 6494.08 Pa and 6499.53 Pa, respectively, with corresponding relative errors of 7.83$\times10^{-4}$ and 1.4$\times10^{-3}$. Therefore, the CSFEM demonstrates higher accuracy than the FEM. In addition, Fig.~\ref{fig:ex02_contour} shows the vertical displacement and stress contours of the CSFEM and FEM at the ultimate bearing capacity. It can be observed that the contour distributions obtained by the CSFEM are smoother than those produced by the FEM.

\begin{table}[H]
  \centering
  \caption{Relative errors in bearing capacity between CSFEM and FEM.}
  \label{tab:ex02:Relative errors}
  \resizebox{\textwidth}{!}{
    \begin{tabular}{ccccc}
      \toprule
      Analytical solution (Pa) & CSFEM (Pa) & FEM (Pa) & Relative error of CSFEM & Relative error of FEM \\
      \midrule
      6489 & 6494.08 & 6499.53 & 7.83$\times10^{-4}$ & 1.4$\times10^{-3}$ \\
      \bottomrule
    \end{tabular}
  }
\end{table}

\subsection{Tunnel excavation}
\begin{figure}[H]
  \centering
  \includegraphics[width=1.0\textwidth]{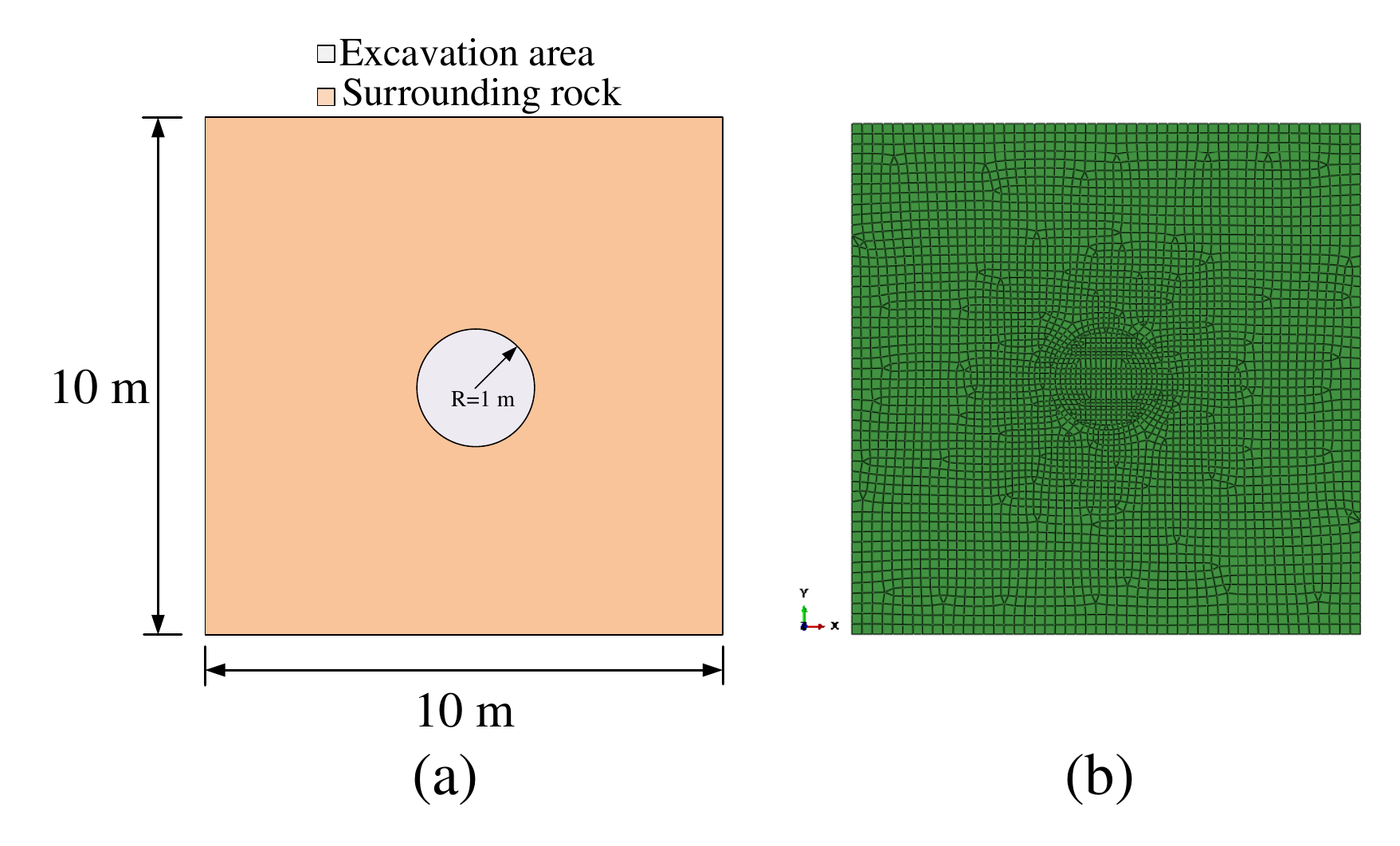}
  \caption{The geometry and mesh model of tunnel excavation; (a) geometry model; (b) mesh model.}
  \label{fig:ex03_geo}
\end{figure}

A tunnel excavation problem is considered in this example, with the geometry shown in Fig.~\ref{fig:ex03_geo}(a). The tunnel has a diameter of 2~m and is excavated in five stages, with each excavation step removing a thickness of 0.4~m. The model is constrained by fixed boundaries at the bottom and normal constraints on the lateral sides. The surrounding rock is characterized by an elastic modulus of $E = 8$~MPa, Poisson’s ratio $\nu = 0.27$, cohesion $c = 120$~kPa, and an internal friction angle of $40^{\circ}$. The analysis begins with an in-situ stress equilibrium step, followed by excavation simulation using the element birth--death technique.

\begin{figure}[H]
  \centering
  \includegraphics[width=0.7\textwidth]{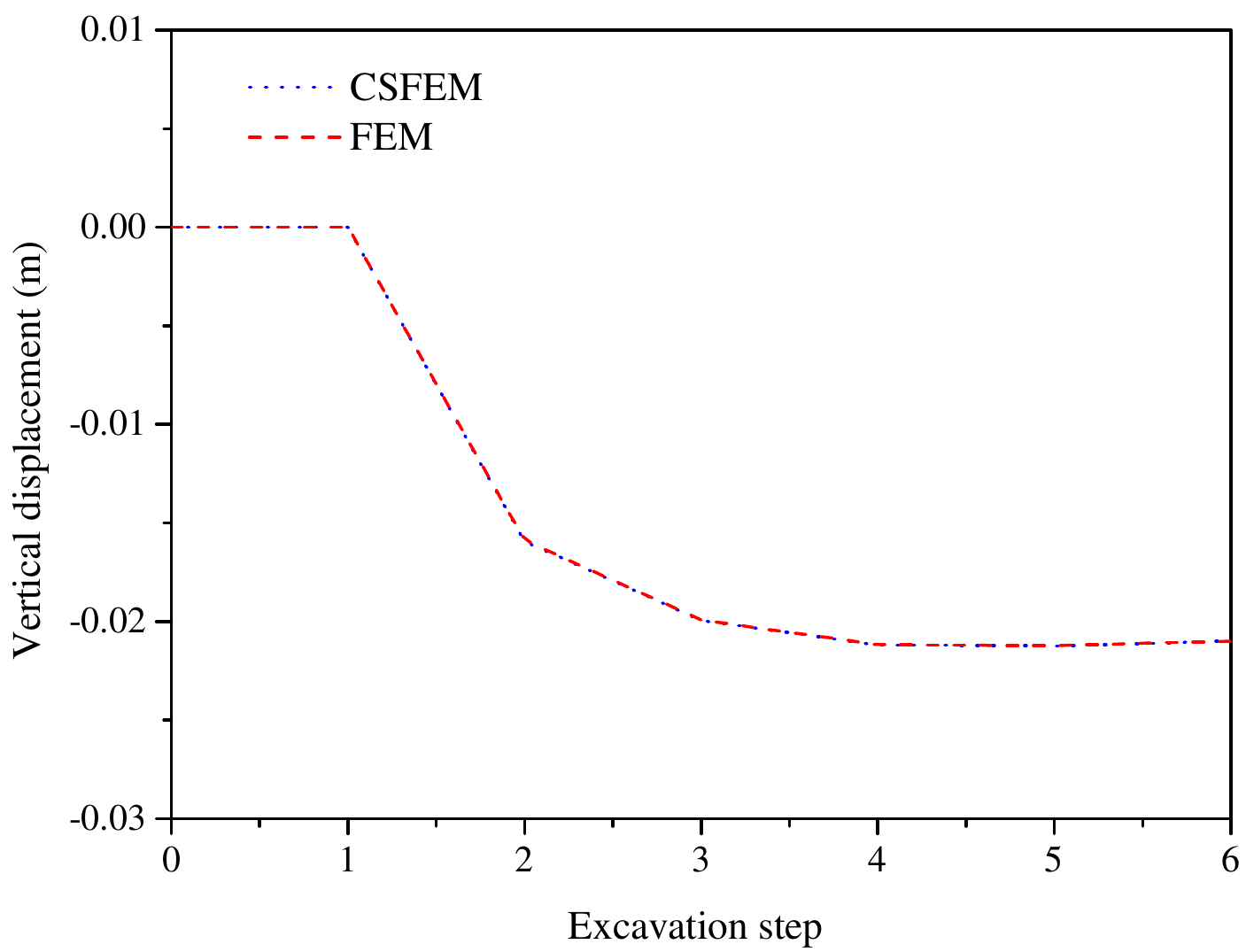}
  \caption{The vertical displacement of CSFEM and FEM under different excavation steps for the crown of the circular tunnel.}
  \label{fig:ex03_U2_kw}
\end{figure}

\begin{figure}[H]
  \centering
  \includegraphics[width=1.0\textwidth]{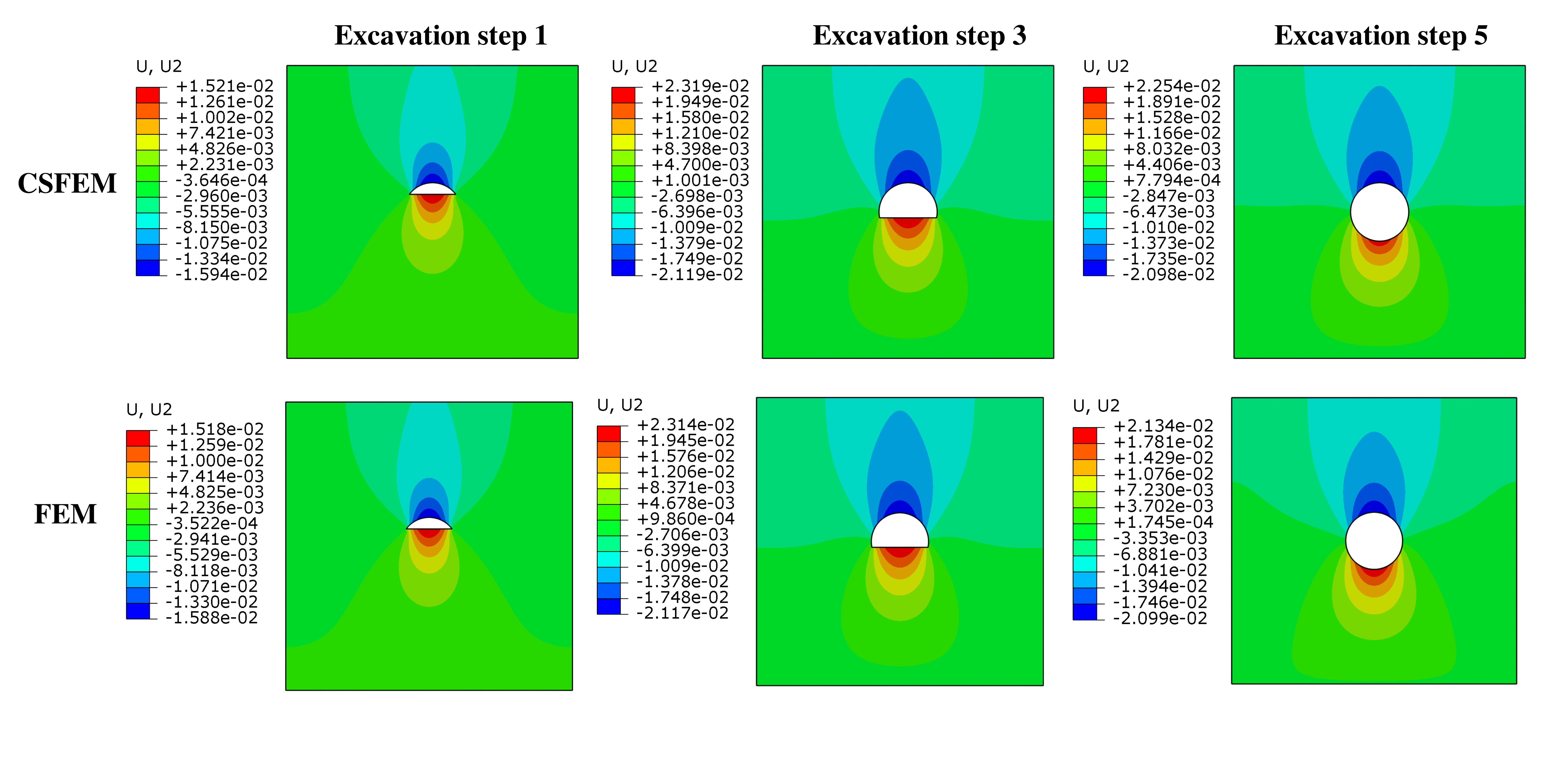}
  \caption{The vertical displacement of CSFEM and FEM under different excavation steps.}
  \label{fig:ex03_contour_U2}
\end{figure}

\begin{figure}[H]
  \centering
  \includegraphics[width=1.0\textwidth]{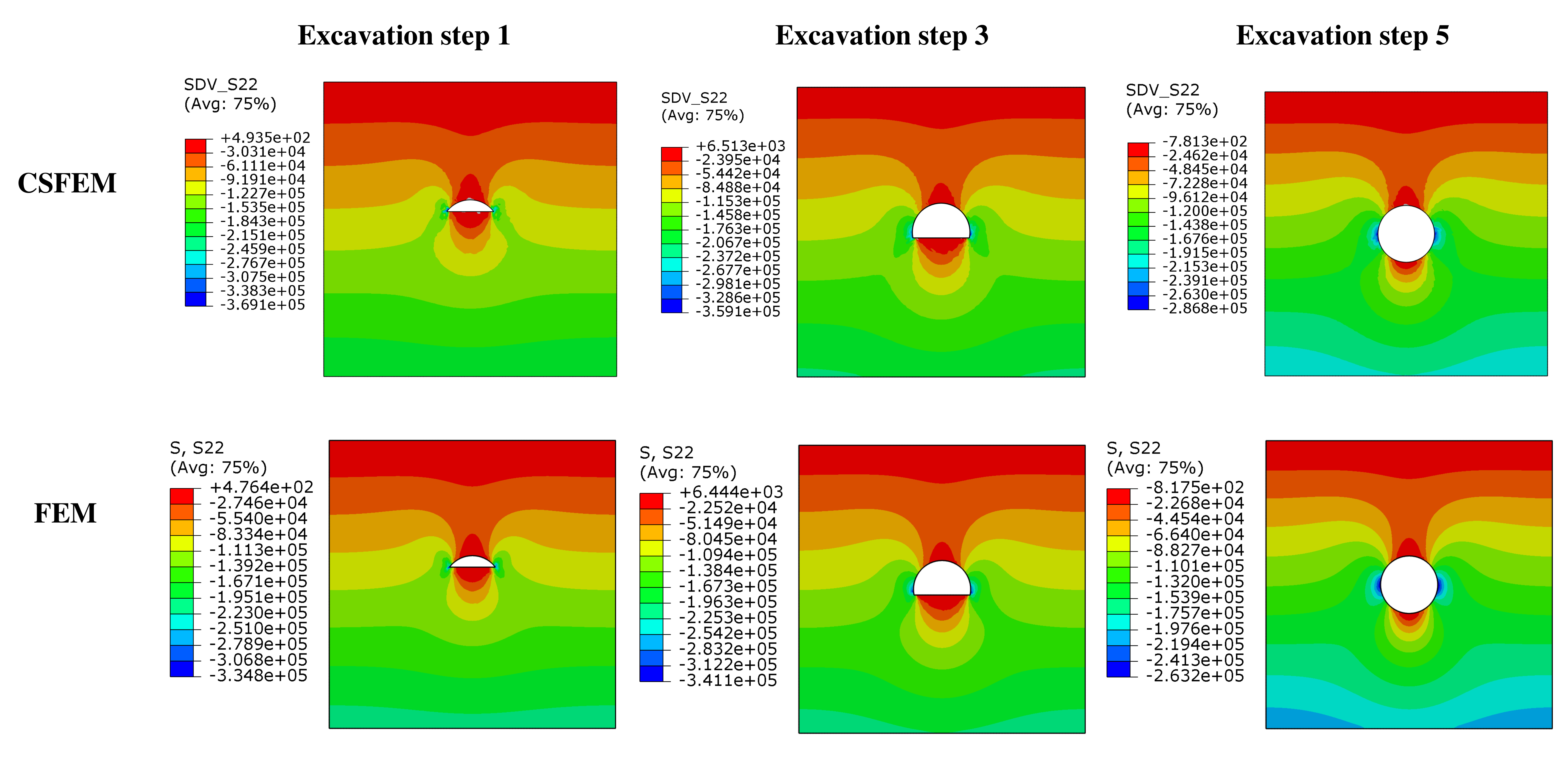}
  \caption{The vertical stress of CSFEM and FEM under different excavation steps.}
  \label{fig:ex03_contour_s22}
\end{figure}

Fig.~\ref{fig:ex03_U2_kw} shows the vertical displacement at the tunnel crown during the different excavation stages, from which it can be observed that CSFEM and FEM exhibit good agreement. When the excavation reaches the fourth stage, the vertical displacement nearly converges, with a maximum value of –0.021 m. Moreover, Fig.~\ref{fig:ex03_contour_U2} and Fig.~\ref{fig:ex03_contour_s22} present the vertical displacement and stress contour plots at the first, third, and fifth excavation stages, further demonstrating the close agreement between the CSFEM and FEM.

\subsection{Slope stability analysis}
\begin{figure}[h]
  \centering
  \includegraphics[width=1.0\textwidth]{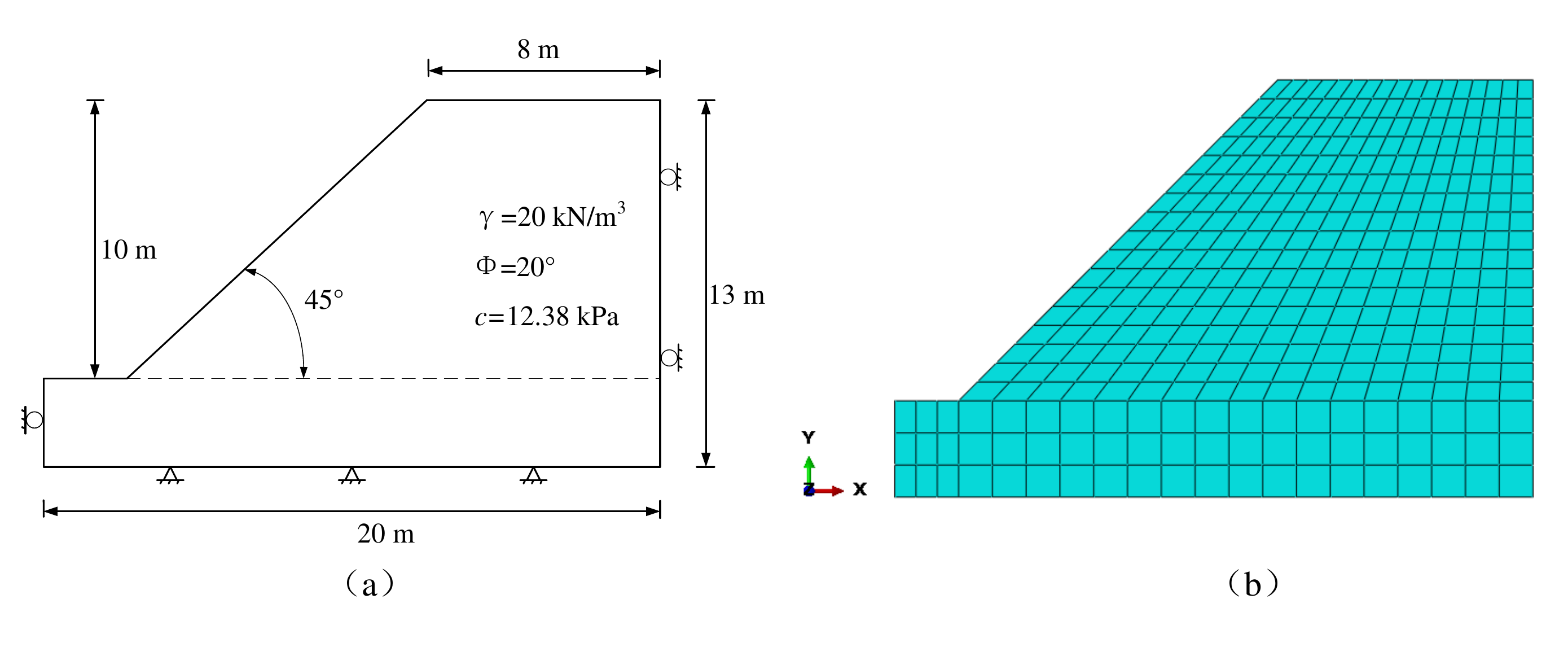}
  \caption{Schematic diagram for the slope stability analysis; (a) geometry and boundary
conditions; (b) mesh model.}
  \label{fig:slope}
\end{figure}

To verify the accuracy of the proposed method, we examined a uniform slope with a height of $H = 10\,\text{m}$, sloping at an angle of $\beta = 45^\circ$ with a friction angle of $\phi = 20^\circ$, a unit weight of $\gamma = 20\,\text{kN/m}^3$, and a cohesion of $c = 12.38\,\text{kPa}$, as shown in Fig.~\ref{fig:slope} (a). Based on the soil properties, the slope has a factor of safety of exactly 1.0, as per the limit analysis solution presented by Chen~\citep{chen2012limit}. The side and bottom boundaries of slopes are constrained without displacement $(\Delta Y = 0)$ and $(\Delta X = 0;\, \Delta Y = 0)$, respectively.

The stability of slope is evaluated using the shear strength reduction technique \citep{matsui1992finite}. The parameters of shear strength are defined as 
\begin{equation}c_m=\frac c{F_r},\end{equation}
\begin{equation}\phi_m=\arctan\biggl(\frac{\tan\phi}{F_r}\biggr),\end{equation}
where $c$ and $\phi$ are the actual shear strength parameters, $c_m$ and $\phi_m$ are the shear strength which can maintain the stability of slope. $F_r$ is the strength reduction factor. In the reduction process, the input shear strength parameters $c$ and $\phi$ are reduction to $c_m$ and $\phi_m$, to trigger the slope failure. The failure criterion was defined by the big jump in the nodal displacement at a chosen point close to the slope’s surface.

\begin{figure}[H]
  \centering
  \includegraphics[width=1.0\textwidth]{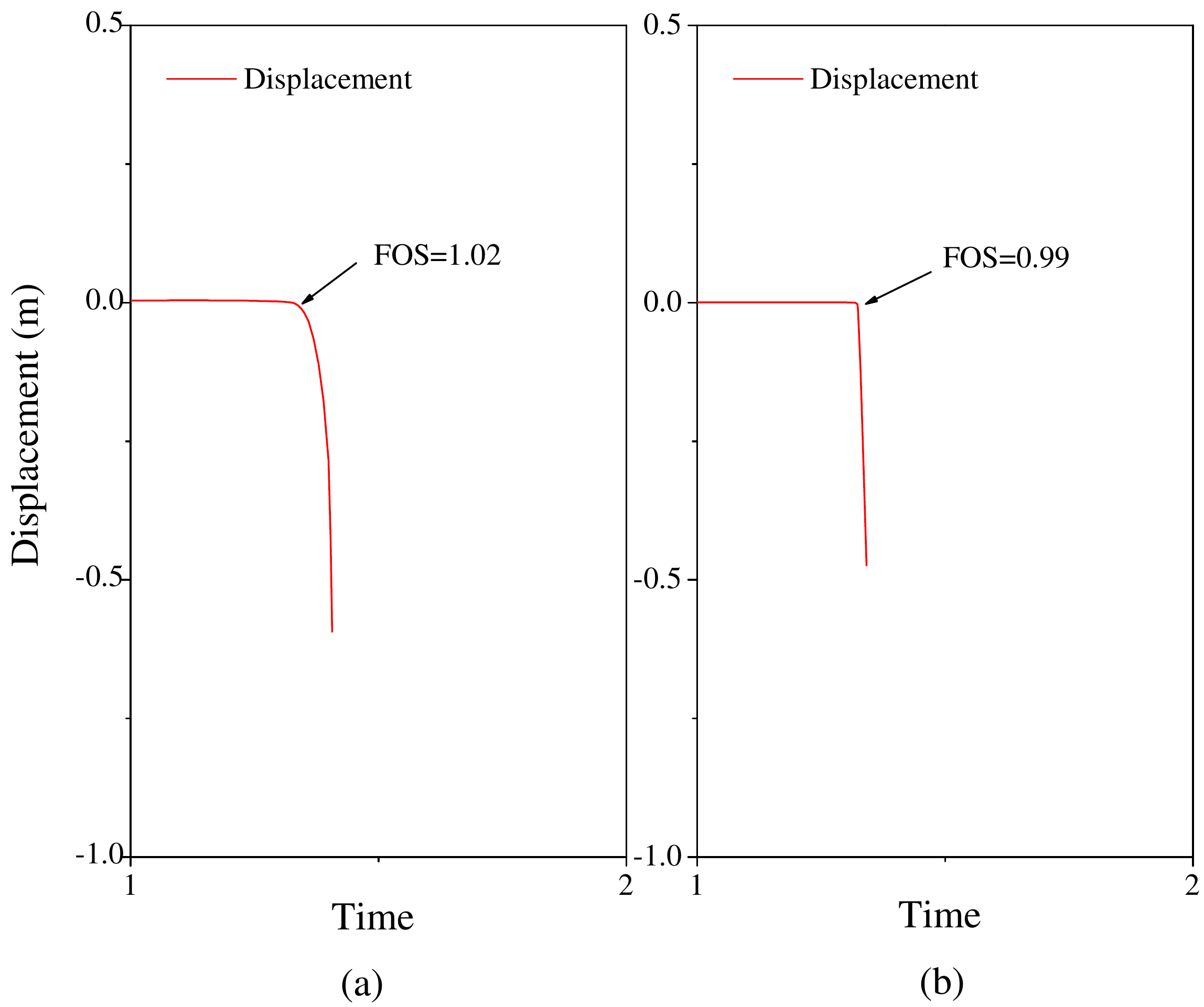}
  \caption{The relationship between displacement and the factor of safety (FOS); (a) the result of CSFEM; (b) the result of FEM.}
  \label{fig:slopejump}
\end{figure}

\begin{figure}[H]
  \centering
  \includegraphics[width=1.0\textwidth]{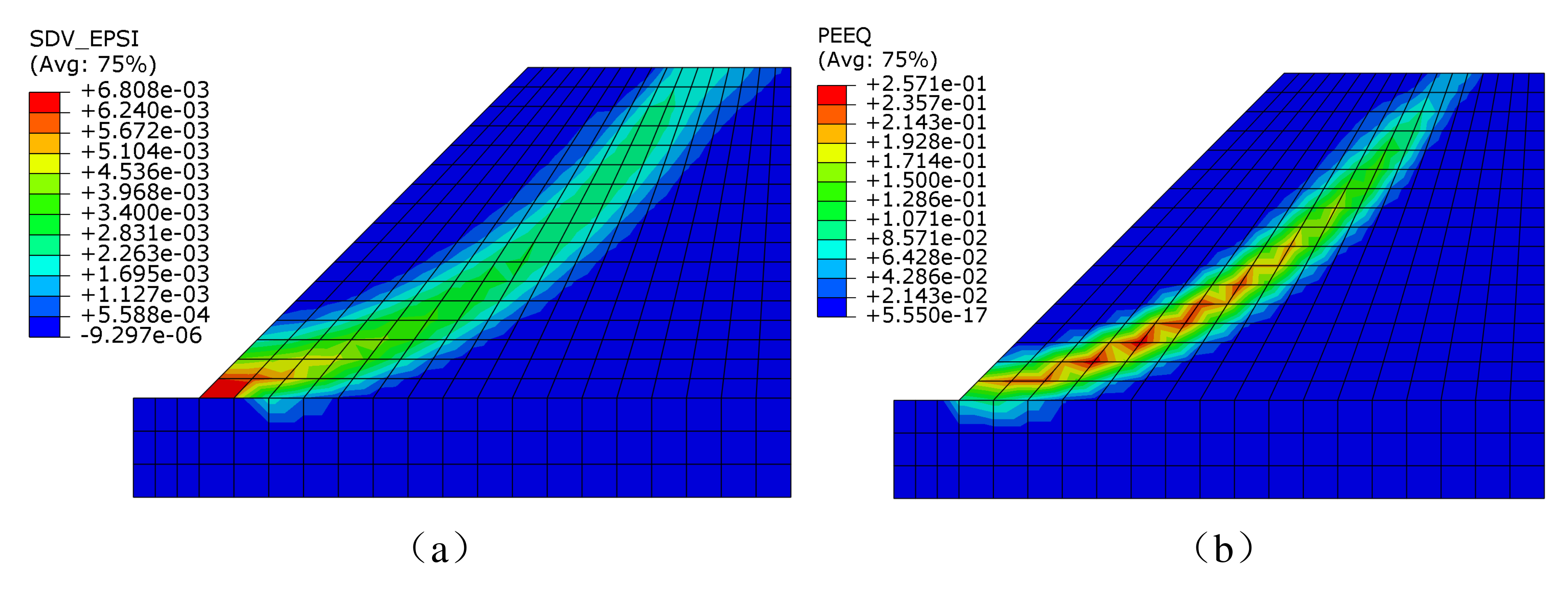}
  \caption{The distribution of plastic strain in the slope at failure
; (a) the result of CSFEM; (b) the result of FEM.}
  \label{ex04:fig:slope_eps}
\end{figure}

We extracted the horizontal displacement of the slope vertex and plotted the relationship between horizontal displacement and the safety factor, as shown in Fig.~\ref{fig:slopejump}. A big jump in displacement occurred in CSFEM when the safety factor reached 1.02, whereas in FEM, the displacement jump occurred at a safety factor of 0.99. It is worth noting that the safety factor obtained by CSFEM is very close to the value of 1.0 reported in the literature ~\citep{dawson1999slope} for this slope. Moreover, Fig.~\ref{ex04:fig:slope_eps} shows the distribution of plastic strain in the slope at failure. It can be seen from the figure that the plastic zones have fully penetrated when the slope fails.

\section{Conclusions}
\label{Conclusions}
An elasto-plastic cell-based smoothed finite element method (CSFEM) has been presented and implemented in ABAQUS for geotechnical analyses. The formulation replaces the compatible strain with a smoothed strain field evaluated over sub-cells and incorporates the Mohr-Coulomb model. The UEL-UMAT data-transfer mechanism is provided to allow direct visualization of stress, strain, and internal variables.

Three benchmark problems, one tunnel excavation, and one slope stability analysis were conducted to assess the performance of the proposed formulation. For the thick-cylinder, biaxial, and strip footing examples, the CSFEM consistently produces smaller displacement or load errors than the conventional FEM, confirming its higher accuracy in both elastic and elasto-plastic regimes. The tunnel excavation and slope stability examples further demonstrate the applicability of the method to staged construction and shear-strength reduction analyses, where the CSFEM and FEM provide comparable global responses, while the CSFEM still delivers well-resolved plastic strain patterns and failure mechanisms consistent with reference solutions.

Overall, the proposed CSFEM offers a simple and effective alternative to conventional FEM for elasto-plastic geotechnical problems. Future work will extend the method to three-dimensional analysis and more advanced constitutive models.

\section{Acknowledgements}
The Yunnan Fundamental Research Projects (grant NO. 202401CF070043), the Xing Dian Talent Support Program of Yunnan Province (XDYC-QNRC-2022-0764) provided support for this study.



 \bibliographystyle{elsarticle-harv} 
\bibliography{cas-refs}





\end{document}